\documentclass[reqno,12pt]{amsart}
\usepackage{amsmath,amsthm,amssymb,amsfonts,amscd}
\usepackage{graphicx,color}

\usepackage{float}
\usepackage{psfrag}
\usepackage{pst-all}

\usepackage{latexsym}

\input xy
\xyoption{all}

\newtheorem{thm}{Theorem}[section]

\newtheorem{cor}[thm]{Corollary}
\newtheorem{lem}[thm]{Lemma}
\newtheorem{prop}[thm]{Proposition}

\theoremstyle{definition}

\newtheorem{defn}[thm]{Definition}%[section]

\newtheorem{rem}[thm]{Remark}

\newtheorem*{pf}{Proof}

\newtheorem{tbl}{Table}%[section]

\numberwithin{equation}{section}
\newcommand{\bp}{\begin{pmatrix}}
\newcommand{\ep}{\end{pmatrix}}
\newcommand{\bps}{\begin{smallmatrix}}
\newcommand{\eps}{\end{smallmatrix}}
\def\C{{\mathbb C}}

\def\Q{{\mathbb Q}}
\def\R{{\mathbb R}}
\def\Z{{\mathbb Z}}
\def\A{{\mathcal A}}
\def\CC{{\mathcal C}}

\def\O{{\mathcal O}}
\def\PP{{\mathcal P}}

\def\Oh{{\hat {\mathcal O}}}
\def\S{{\mathcal S}}

\def\ZZ{{\mathcal Z}}

\def \fC{\frak C}
\def \fS{{\frak S}}

\def\p{{\partial }}

\def\i{\sqrt{-1}}

\def\bb{{\bar b}}

\def\sb{{\bar s}}

\def\qb{{\bar q}}

\def \fpartial#1{\frac{\partial}{\partial #1}}

\def \0{{\bf 0}}
\def \1{{\bf 1}}
\def \bn{{\bf n}}

\def\Tr{\mathrm{Tr}}
\def\Hom{\mathrm{Hom}}
\def \Irr{\mathrm{Irr}}

\def \Mat{\mathrm{Mat}}

\def \Ob{\mathrm{Ob}}
\def \diag{\mathrm{diag}}

\def \Id{\mathrm{id}}
\def \rad{\mathrm{rad}}
\def \Irr{\mathrm{Irr}}
\def \mod{\mathrm{mod}\,\text{-}\,}

\def \lrDelta{\overset{\leftrightarrow}{\Delta}}

\def \mf#1#2#3#4{
\xymatrix{{#1}\  \ar@<0.4ex>[r]^{{#2}} & \ {#4}
\ar@<0.4ex>[l]^{{#3}}
}
}

\def \mfs#1#2#3#4{\!
\xymatrix@C=1,5em{{#1} \! \ar@<0.2ex>[r]^{{#2}} & \! {#4}
\ar@<0.2ex>[l]^{{#3}}
}
\!}

\def \mfl#1#2#3#4{
\xymatrix@C=2.6em{{#1}\  \ar@<0.4ex>[r]^{{#2}} &\  {#4}
\ar@<0.2ex>[l]^{{#3}}
}
}

\def \mfss#1#2#3#4{\!
\xymatrix@C=1.5em{{#1} \ar@<0.3ex>[r]^{{#2}} & {#4}
\ar@<0.3ex>[l]^{{#3}}
}
\!}

\def \Abn{{\C\vec{\Gamma}(\bn)}}

\def \MF#1{\mathit{MF_{#1}(f)}}
\def \tMF#1{\mathit{MF^{gr}_{#1}(f)}}
\def \HMF#1{HMF_{#1}(f)}
\def \tHMF#1{{HMF}^{gr}_{#1}(f)}

\def \ti#1{\widetilde{#1}}

\def \HomD#1{\Hom_{\HMF{#1}}}
\def \HomDz#1{\Hom_{\tHMF{#1}}}

\def \ov#1{\frac{1}{#1}}

\def \ie{{\it i.e.}}

\newcommand{\qgr}{\operatorname{\mathrm{qgr}}}
\newcommand{\vecx}{\ensuremath{\vec{x}}}
\newcommand{\vx}{\vecx}
\newcommand{\vecc}{\ensuremath{\vec{c}}}
\newcommand{\vc}{\vecc}
\newcommand{\gr}{\ensuremath{\mathrm{gr\text{--}}}}
\newcommand{\tor}{\ensuremath{\mathrm{tor\text{--}}}}
\newcommand{\scO}{\ensuremath{\mathcal{O}}}
\newcommand{\vecomega}{\ensuremath{\vec{\omega}}}
\def \bZ{\Z}
\newcommand{\Dbsing}{\ensuremath{D^{\mathrm{gr}}_{\mathrm{Sg}}}}
\newcommand{\grproj}{\ensuremath{\mathrm{grproj\text{--}}}}
\newcommand{\End}{\mathop{\mathrm{End}}\nolimits}

\newcommand{\pp}{\ensuremath{\boldsymbol{p}}}

\begin{document}
\begin{flushleft}
\hfill YITP-05-65\\
\hfill RIMS-1521\, \\
\end{flushleft}

\vspace*{0.4cm}

\title{Matrix Factorizations and Representations of Quivers II:\ 
type ADE case}

\date{November 7th, 2005}
\author{Hiroshige Kajiura}
\address{Yukawa Institute for Theoretical Physcics, Kyoto University,
Kyoto 606-8502, Japan}
\email{kajiura@yukawa.kyoto-u.ac.jp}
\author{Kyoji Saito}
\address{Research Institute for Mathematical Sciences, Kyoto University,
Kyoto 606-8502, Japan}
\email{saito@kurims.kyoto-u.ac.jp}
\author{Atsushi Takahashi}
\address{Research Institute for Mathematical Sciences, Kyoto University,
Kyoto 606-8502, Japan}
\email{atsushi@kurims.kyoto-u.ac.jp}
\begin{abstract}
We study a triangulated category of graded matrix factorizations 
for a polynomial of type ADE. 
We show that it is equivalent to the derived category of finitely
generated modules over the path algebra of the corresponding Dynkin quiver. 
Also, we discuss a special stability condition for the triangulated 
category in the sense of T.~Bridgeland, which is naturally defined 
by the grading. 
\end{abstract}

\maketitle
%
%
%%%%%%%%%%%%%%%%%%%%%%%%%%%%%%%%%%%%%%%%%%%%%%%%%%%%%%%%%%%%%%%%%%%%%%%%%%%%%%
%%%%%%%%%%%%%%%%%%%%%%%%%%%%%%%%%%%%%%%%%%%%%%%%%%%%%%%%%%%%%%%%%%%%%%%%%%%%%%
%%%%%%%%%%%%%%%%%%%%%%%%%%%%%%%%%%%%%%%%%%%%%%%%%%%%%%%%%%%%%%%%%%%%%%%%%%%%%%

\tableofcontents

\section{Introduction}

The universal deformation and the simultaneous resolution of a simple 
singularity are described by the corresponding simple Lie algebra 
(Brieskorn \cite{bs:1}). 
Inspired by that theory, the second named author 
associated in \cite{sa:elliptic},\cite{sa:around} 
a generalization of root systems, 
consisting of vanishing cycles of the singularity, 
to any regular weight systems \cite{sa:weight}, 
and asked to construct a suitable Lie theory 
in order to reconstruct the primitive forms 
for the singularities. 
In fact, the simple singularities correspond exactly to 
the weight systems having only positive exponents, 
and, in this case, 
this approach gives the classical finite root systems 
as in \cite{bs:1}. 
As the next case, 
the approach is worked out for simple 
elliptic singularities corresponding to weight 
systems having only non-negative exponents, 
from where the theory of elliptic Lie algebras 
is emerging \cite{sa:elliptic}. 
However, the root system in this approach 
in general is hard to manipulate 
because of the transcendental nature of vanishing cycles.
Hence, 
he asked (\cite{sa:around}, Problem in p.124 in English version) 
{\em an algebraic and/or a combinatorial 
construction of the root system starting from a regular weight system}.

In \cite{t:2}, 
based on the mirror symmetry for the Landau-Ginzburg orbifolds 
and also based on the duality theory of the weight systems 
\cite{sa:duality},\cite{t:1}, 
the third named author proposed a new approach to the root systems, 
answering to the above problem. 
He introduced a triangulated category $D^b_\Z(\A_f)$ 
of graded matrix factorizations 
for a weighted homogeneous polynomial $f$ 
attached to a regular weight system 
and showed that the category $D^b_\Z(\A_f)$ for 
a polynomial of type $A_l$ 
is equivalent to the bounded derived 
category of modules over the path algebra of the Dynkin quiver 
of type $A_l$. 
He conjectured (\cite{t:2}, Conjecture 1.3) 
further that the same type of equivalences hold 
for all simple polynomials 
of type ADE. The main goal of the present paper is to answer affirmatively 
to the conjecture. 

One side of this conjecture: 
the properties of the category of 
modules over a path algebra of a Dynkin quiver 
are already well-understood 
by the Gabriel's theorem \cite{g:1}, 
which states that the number of the indecomposable objects 
in the category for a Dynkin quiver 
coincides with the number of the positive roots of the root system 
corresponding to the Dynkin diagram. 
The other side of the conjecture: 
the triangulated categories of (ungraded) 
matrix factorizations were introduced and 
developed by Eisenbud \cite{e:1} and Kn\"orrer \cite{k:1} 
in the study of the maximal Cohen-Macaulay modules. 
Recently, the categories of matrix factorizations are rediscovered 
in string theory as the categories of topological D-branes of type B 
in Landau-Ginzburg models (see \cite{kl:1},\cite{kl:2}). 
The category $D^b_\Z(\A_f)$ of graded matrix factorizations 
is then motivated 
by the work on the categories of topological D-branes of type B 
in Landau-Ginzburg {\em orbifolds} $(f,\Z/h\Z)$ by Hori-Walcher
\cite{hw:1}, where the orbifolding corresponds to introducing the 
$\Q$-grading. 
In fact, in \cite{t:2}, 
the triangulated category $D^b_\Z(\A_f)$ is constructed 
from a special $A_\infty$-category with $\Q$-grading 
via the twisted complexes 
in the sense of Bondal-Kapranov \cite{bk:1}. 
Independently, D.~Orlov defines a triangulated category, 
called the category of graded D-branes of type B, 
which is in fact equivalent to $D^b_\Z(\A_f)$ 
(see the end of subsection \ref{ssec:Dz}). 
Though some notions, for instance the central charge of the stability
condition (see section \ref{sec:stab}), 
can be understood more naturally in $D^b_\Z(\A_f)$, 
the Orlov's construction of categories 
requires less terminologies and is easier to understand 
in a traditional way in algebraic geometry. 
Therefore, in this paper 
we shall use the Orlov's construction with a slight modification 
of the scaling of degrees and denote 
the modified category by $\tHMF{R}$.

Let us explain details of the contents of the present paper. 
In section \ref{sec:cat}, we recall the construction of
triangulated categories of matrix factorizations. 
Since we compare the category $HMF_R^{gr}(f)$ 
with the ungraded version $HMF_{\mathcal O}(f)$ 
in the proof of our main theorem (Theorem \ref{thm:main}), 
we first introduce the ungraded version 
$HMF_{\mathcal O}(f)$ corresponding to that given in \cite{o:1} 
in subsection \ref{ssec:D}, and then 
we define the graded version $HMF_R^{gr}(f)$ 
based on \cite{o:2} in subsection \ref{ssec:Dz}, where 
we also explain the relation of the category $HMF_R^{gr}(f)$ with 
the category $D^b_\Z(\A_f)$ introduced in \cite{t:2}. 
Section \ref{sec:main} is the main part of the present paper. 
In subsection \ref{ssec:main-state}, we state the main theorem 
(Theorem \ref{thm:main}): 
for a polynomial $f$ of type ADE, 
$HMF_R^{gr}(f)$ is equivalent as a triangulated category 
to the bounded derived category of modules over the path algebra 
of the Dynkin quiver of type of $f$. 
Subsection \ref{ssec:proof} is devoted to the proof of 
Theorem \ref{thm:main}. 
The proof is based on various explicit data 
on the matrix factorizations; 
the complete list of the matrix factorizations 
(Table \ref{tbl:1}), their gradings (Table \ref{tbl:2}) and 
the complete list of the morphisms in $HMF_R^{gr}(f)$ 
(Table \ref{tbl:3}). 
The tables are arranged in the final section (section \ref{sec:tbl}). 
In section \ref{sec:stab}, 
we construct a stability condition, 
the notion of which is introduced by Bridgeland \cite{bd:1}, 
for the triangulated category $HMF_R^{gr}(f)$. 
One can see that, as in the $A_l$ case \cite{t:2}, 
the phase of objects 
(see Theorem \ref{thm:grading} or Table \ref{tbl:2}) 
and the central charge $\ZZ$ (Definition \ref{defn:centralcharge}) 
can be naturally given by the grading of matrix factorizations 
in Table \ref{tbl:2} (c.f. \cite{w:1}). 
They in fact define a stability condition 
on $HMF_R^{gr}(f)$ (Theorem \ref{thm:stab}), 
from which an abelian category is obtained as 
a full subcategory of $HMF_R^{gr}(f)$. 
In Proposition \ref{prop:abelian}, we show 
that this abelian category is equivalent to an abelian category of 
modules over the path algebra $\C\vec{\Delta}_{principal}$, 
where $\vec{\Delta}_{principal}$ is 
the Dynkin quiver with the orientation being taken to be 
the principal orientation introduced in \cite{sa:principal}. 
In the appendix, we include another proof of Theorem \ref{thm:main} 
by K.~Ueda. 

{\bf Acknowledgement} :\ 
We are grateful to M.~Kashiwara and K.~Watanabe for valuable discussions.
This work was partly supported by Grant-in Aid for Scientific Research 
grant numbers 16340016, 17654015 and 17740036 of the Ministry of Education, 
Science and Culture in Japan. 
H.~K is supported by JSPS Research Fellowships for Young 
Scientists.

%%%%%%%%%%%%%%%%%%%%%%%%%%%%%%%%%%%%%%%%%%%%%%%%%%%%%%%%%%%%%%%%%%%%%%%%%%%%%%

\section{Triangulated categories of matrix factorizations}
\label{sec:cat}

In this section, we set up several definitions which are used 
in the present paper. 
The goal of this section is the introduction of the categories 
$\HMF{A}$ and $\tHMF{A}$ attached 
to a weighted homogeneous polynomial $f\in A$ 
following \cite{o:1},\cite{o:2} with slight modifications.

 \subsection{The triangulated category ${\HMF{A}}$ 
of matrix factorizations} \hfill \\
\label{ssec:D}
Let $A$ be either the polynomial ring $R:=\C[x,y,z]$, 
the convergent power series ring $\O:=\C\{x,y,z\}$ or 
the formal power series ring $\Oh:=\C[[x,y,z]]$ 
in three variables $x,y$ and $z$. 
\begin{defn}[Matrix factorization]
For a nonzero polynomial $f\in A$, 
a {\em matrix factorization} $M$ of $f$ is defined by 
\begin{equation*}
 M:=\Big( \mf{P_0}{p_0}{p_1}{P_1}\Big) \ ,
\end{equation*}
where $P_0$, $P_1$ are right {\em free} 
$A$-modules of finite rank, 
and $p_0:P_0\to P_1$, $p_1:P_1\to P_0$ are 
$A$-homomorphisms such that $p_1p_0=f\cdot\Id_{P_0}$ and
$p_0p_1=f\cdot\Id_{P_1}$. 
The set of all matrix factorizations of $f$ 
is denoted by $\MF{A}$. 
 \label{defn:MF}
\end{defn}
Since $p_0p_1$ and $p_1p_0$ are 
$f$ times 
the identities, 
where $f$ is nonzero element of $A$, 
the rank of $P_0$ coincides with that of $P_1$. 
We call the rank the {\em size} of the matrix factorization $M$. 
\begin{defn}[Homomorphism]
Given two matrix factorizations 
$M:=(\mfs{P_0}{p_0}{p_1}{P_1})$ and 
$M':=(\mfs{P'_0}{p'_0}{p'_1}{P'_1})$, 
a {\em homomorphism} $\Phi:M\to M'$ 
is a pair of $A$-homomorphisms $\Phi=(\phi_0,\phi_1)$ 
\begin{equation*}
 \phi_0:P_0\to P_0'\ ,\qquad \phi_1:P_1\to P_1'\ ,
\end{equation*}
such that the following diagram commutes: 
\begin{equation*}
 \xymatrix{
 \ P_0\ \ar[r]^{p_0} \ar[d]^{\phi_0} & 
 \ P_1\ \ar[r]^{p_1} \ar[d]^{\phi_1} & 
 \ P_0\ \ar[d]^{\phi_0} & \\
 \ P'_0\ \ar[r]^{p'_0} & 
 \ P'_1\ \ar[r]^{p'_1} & 
 \ P'_0\ & \hspace*{-1.7cm} .
}
\end{equation*}
The set of all homomorphisms from $M$ to $M'$, 
denoted by $\Hom_\MF{A}(M,M')$, is naturally an $A$-module and 
is finitely generated, 
since the sizes of the matrix factorizations are finite. 
For three matrix factorizations $M$, $M'$, $M''$ 
and homomorphisms $\Phi:M\to M'$ and $\Phi':M'\to M''$, 
the composition $\Phi'\Phi$ is defined by 
\begin{equation*}
 \Phi'\Phi=(\phi'_0\phi_0,\phi_1'\phi_1)\ .
\end{equation*}
This composition 
is associative: $\Phi''(\Phi'\Phi)=(\Phi''\Phi')\Phi$ 
for any three homomorphisms. 
 \label{defn:hom}
\end{defn}
\begin{defn}[$\HMF{A}$]
An additive category $\HMF{A}$ is defined by the following data. 
The set of objects is given by the set of all matrix factorizations: 
\begin{equation*}
 \Ob(\HMF{A}):=\MF{A}\ . 
\end{equation*}
For any two objects $M,M'\in\MF{A}$, 
the set of morphisms is given by the quotient module: 
\begin{equation*}
 \HomD{A}(M,M'):=
 \Hom_\MF{A}(M,M')/ \sim \ ,
\end{equation*}
where two elements $\Phi$, $\Phi'$ in $\Hom_\MF{A}(M,M')$ are equivalent 
(homotopic) 
$\Phi\sim\Phi'$ if there exists 
a {\em homotopy} $(h_0,h_1)$, \ie, 
a pair $(h_0,h_1):(P_0\to P_1',P_1\to P_0')$ of $A$-homomorphisms 
such that $\Phi'-\Phi=(p_1'h_0+h_1p_0, p_0'h_1+h_0p_1)$. 
The composition of morphisms on $\HomD{A}$ is induced from that 
on $\Hom_\MF{A}$ 
since $\Phi\sim\Phi'$ and $\Psi\sim\Psi'$ imply $\Psi\Phi\sim\Psi'\Phi'$. 
\end{defn}
Note that the matrix factorization $M=(\mfs{P_0}{p_0}{p_1}{P_1})\in\MF{A}$ 
of size one with $(p_0,p_1)=(1,f)$ or $(p_0,p_1)=(f,1)$ 
defines the zero object in $\HMF{A}$, 
that is: one has 
$\HomD{A}(M,M')=\HomD{A}(M',M)=0$ 
for any matrix factorization $M'\in\MF{A}$ 
or, equivalently, $\Id_M\in\Hom_{\MF{A}}(M,M)$ is homotopic to zero. 
\begin{lem}
For any two matrix factorizations $M,M'\in\HMF{A}$, 
the space of morphisms $\HomD{A}(M,M')$ is a finitely generated 
$A/\left(\frac{\p f}{\p x},\frac{\p f}{\p y},\frac{\p f}{\p z}\right)$-module. 
 \label{lem:finite}
\end{lem}
\begin{pf}
Since 
$\Hom_\MF{A}(M,M')$ is a finitely generated $A$-module and 
the equivalence relation $\sim$ is given by 
quotienting out by an $A$-submodule, 
$\HomD{A}(M,M')$ is also a finitely generated $A$-module. 
On the other hand, the Jacobi ideal 
$\left(\frac{\p f}{\p x},\frac{\p f}{\p y},\frac{\p f}{\p z}\right)$ 
annihilates $\HomD{A}(M,M')$: that is, 
$\frac{\p f}{\p x}\left(\phi_0,\phi_1\right)\sim 
\frac{\p f}{\p y}\left(\phi_0,\phi_1\right)\sim 
\frac{\p f}{\p z}\left(\phi_0,\phi_1\right)\sim 0$ 
for any morphism $\Phi=(\phi_0,\phi_1)$. 
This can be shown for instance by 
differentiating $p_1p_0=f\cdot\Id_{P_0}$ and 
$p_0p_1=f\cdot\Id_{P_1}$ 
by $\fpartial{x}$. Then 
we have two identities 
$\frac{\p p_1}{\p x}p_0+p_1\frac{\p p_0}{\p x}=
\frac{\p f}{\p x}\cdot\Id_{P_0}$ 
and 
$\frac{\p p_0}{\p x}p_1+p_0\frac{\p p_1}{\p x}
=\frac{\p f}{\p x}\cdot\Id_{P_1}$. 
Multiplying $\phi_0$ and $\phi_1$ by these two identities, respectively, 
leads to $\frac{\p f}{\p x}\left(\phi_0,\phi_1\right)\sim 0$, 
where $(\phi_1\frac{\p p_0}{\p x},\phi_0\frac{\p p_1}{\p x})$ 
is the corresponding homotopy. 
In a similar way one can obtain 
$\frac{\p f}{\p y}\left(\phi_0,\phi_1\right)\sim 
\frac{\p f}{\p z}\left(\phi_0,\phi_1\right)\sim 0$. 
\qed\end{pf}
\begin{defn}[Shift functor]
The {\em shift functor} $T:\HMF{A}\to\HMF{A}$ 
is defined as follows. 
The action of $T$ on 
$M=(\mfs{P_0}{p_0}{p_1}{P_1})\in\HMF{A}$ is given by 
\begin{equation*}
 T\Big( \mf{P_0}{p_0}{p_1}{P_1}\Big)
 :=\Big(\mf{P_1}{-p_1}{-p_0}{P_0}\Big)\ .
\end{equation*}
For any $M,M'\in\HMF{A}$, the action of $T$ on 
$\Phi=(\phi_0,\phi_1)\in\HomD{A}(M,M')$ is given by 
\begin{equation*}
 T(\phi_0,\phi_1):=(\phi_1,\phi_0)\ . 
\end{equation*}
 \label{defn:shift}
\end{defn}
Note that 
the square $T^2$ of the shift functor is isomorphic to the identity functor 
on $\HMF{A}$. 
\begin{defn}[Mapping cone]
For an element $\Phi=(\phi_0,\phi_1)\in\Hom_\MF{A}(M,M')$, 
the {\em mapping cone} $C(\Phi)\in\MF{A}$ is defined by 
\begin{equation*}
 \begin{split}
 & C(\Phi):=\Big( \mf{C_0}{c_0}{c_1}{C_1}\Big)\ ,\quad\text{where} \\
 & C_0:=P_1\oplus P_0'\ ,\qquad C_1:=P_0\oplus P_1'\ ,\qquad 
   c_0:=\bp -p_1 & 0 \\ \phi_1 & p_0'\ep\ ,\qquad 
   c_1:=\bp -p_0 & 0 \\ \phi_0 & p_1'\ep\ .
 \end{split}
\end{equation*}
 \label{defn:cone}
\end{defn}
The following is stated in \cite{o:1} Proposition 3.3. 
\begin{prop}
The additive category $\HMF{A}$ endowed with the shift functor 
$T$ and the distinguished triangles 
forms a triangulated category, 
where a distinguished triangle is 
a sequence of morphisms which is isomorphic to the 
sequence 
\begin{equation*}
 M\stackrel{\Phi}{\to} M' \to C(\Phi)\to T(M)\ 
\end{equation*}
for some $M,M'\in\MF{A}$ and $\Phi\in\Hom_\MF{A}(M,M')$. 
\end{prop}
\begin{pf}
The proof is 
the same as the proof of the analogous 
result for a usual homotopic category 
(see e.g. \cite{gm:1}, \cite{ks:1}). 
\qed\end{pf}

 \subsection{The triangulated category $\tHMF{R}$ 
of graded matrix factorizations}
\hfill\\
\label{ssec:Dz}
In this subsection, 
we study graded matrix factorizations for 
a weighted homogeneous polynomial $f$ 
and construct the corresponding triangulated category, denoted by 
$\tHMF{R}$. 

A quadruple $W:=(a,b,c;h)$ of positive integers 
with $\mathrm{g.c.d}(a,b,c)=1$ is called 
a {\em weight system}. 
For a weight system $W$, we define the 
{\em Euler vector field} $E=E_W$ by 
\begin{equation*}
  E:= \frac{a}{h}x\frac{\p}{\p x}+\frac{b}{h}y\frac{\p}{\p y}+
 \frac{c}{h}z\frac{\p}{\p z}\ .
\end{equation*}
For a given weight system $W$, 
${R}$ becomes a graded ring by putting 
$\deg(x)=\frac{2a}{h}$, $\deg(y)=\frac{2b}{h}$ and 
$\deg(z)=\frac{2c}{h}$. 
Let ${R}=\oplus_{d\in\frac{2}{h}\Z_{\ge 0}} {R}_d$ 
be the graded piece decomposition, 
where 
${R}_d:=\{f\in R\ |\ 2Ef=df\}$\ .
A weight system $W$ is called {\em regular} (\cite{sa:weight}) 
if the following equivalent conditions are satisfied: 
\begin{itemize}
 \item[(a)] $\chi_W(T):=
T^{-h}\frac{(T^h-T^a)(T^h-T^b)(T^h-T^c)}{(T^a-1)(T^b-1)(T^c-1)}$ 
has no poles except at $T=0$. 

 \item[(b)] 
A generic element of the eigenspace $R_2=\{f\in R\ |\ Ef= f\}$ 
has an isolated critical point at the origin, 
{\it i.e.}, 
the Jacobi ring 
$R/\left(\frac{\p f}{\p x},\frac{\p f}{\p y},\frac{\p f}{\p z}\right)$ 
is finite dimensional over $\C$. 
\end{itemize}
Such an element $f$ of $R_2$ as in (b) 
shall be called a {\em polynomial of type $W$}.

In the present paper, 
{\em by a graded module, we mean a graded right module with degrees 
only in $\frac{2}{h}\Z$}. 
Namely, a graded ${R}$-module $\ti{P}$ decomposes 
into the direct sum: 
\begin{equation}
\ti{P}=\oplus_{d\in\frac{2}{h}\Z}\ti{P}_d\ .
 \label{P-grad}
\end{equation}
For two graded ${R}$-modules $\ti{P}$ and $\ti{P}'$, 
a graded ${R}$-homomorphism $\phi$ of degree $s\in\frac{2}{h}\Z$ 
is an ${R}$-homomorphism $\phi:\ti{P}\to\ti{P}'$ such that 
$\phi(\ti{P}_d)\subset\ti{P}'_{d+s}$ for any $d$. 
The category of graded ${R}$-modules has a degree shifting 
automorphism $\tau$ defined by
\footnote{This $\tau$ is what is often denoted 
(for instance \cite{y},\cite{o:2}) by $(1)$, 
{\it i.e.} , $\tau(\ti{P})=\ti{P}(1)$. 
} 
\begin{equation*}
 (\tau(\ti{P}))_d:=\ti{P}_{d+\frac{2}{h}}\ .
\end{equation*}
For any two graded ${R}$-modules $\ti{P}$, $\ti{P}'$ 
and a graded ${R}$-homomorphism $\phi:\ti{P}\to\ti{P}'$, 
we denote the induced graded ${R}$-homomorphism by 
$\tau(\phi):\tau(\ti{P})\to\tau(\ti{P}')$. 
On the other hand, an ${R}$-homomorphism $\phi:\ti{P}\to\ti{P}'$ 
of degree $\frac{2m}{h}$ induces 
a degree zero ${R}$-homomorphism from $\ti{P}$ to $\tau^m(\ti{P}')$, 
which we denote again by $\phi:\ti{P}\to\tau^m(\ti{P}')$. 
\begin{defn}[Graded matrix factorization]
For a 
polynomial $f$ of type $W$, 
a {\em graded matrix factorization} $\ti{M}$ of $f\in {R}$ 
is defined by 
\begin{equation*}
 \ti{M}:=\Big( \mf{\ti{P}_0}{p_0}{p_1}{\ti{P}_1}\Big) \ ,
\end{equation*}
where $\ti{P}_0$, $\ti{P}_1$ are 
free graded right ${R}$-modules of finite rank, 
$p_0:\ti{P}_0\to\ti{P}_1$ is a graded ${R}$-homomorphism of degree zero, 
$p_1:\ti{P}_1\to\ti{P}_0$ is a graded ${R}$-homomorphism of degree two 
such that $p_1p_0=f\cdot\Id_{\ti{P}_0}$ and $p_0p_1=f\cdot\Id_{\ti{P}_1}$. 
The set of all graded matrix factorizations of $f$ 
is denoted by $\tMF{R}$. 
 \label{defn:gMF}
\end{defn}
\begin{defn}[Homomorphism]
Given two graded matrix factorizations 
$\ti{M},\ti{M}'\in\tMF{R}$, 
a {\em homomorphism} 
$\Phi=(\phi_0,\phi_1):\ti{M}\to\ti{M}'$ 
is a homomorphism 
in the sense of Definition \ref{defn:hom} 
such that $\phi_0$ and $\phi_1$ are graded ${R}$-homomorphisms 
of {\em degree zero}. 
The vector space of all graded 
${R}$-homomorphisms from $\ti{M}$ to $\ti{M}'$ 
is denoted by $\Hom_\tMF{R}(\ti{M},\ti{M}')$. 
 \label{defn:gHom}
\end{defn}
For three graded matrix factorizations 
$\ti{M},\ti{M}',\ti{M}''\in\tMF{R}$ and 
morphisms $\Phi:\ti{M}\to\ti{M}'$, 
$\Phi':\ti{M}'\to\ti{M}''$, the composition is again 
a graded ${R}$-homomorphism: 
$\Phi'\Phi\in\Hom_\tMF{R}(\ti{M},\ti{M}'')$. 
\begin{defn}[$\tHMF{R}$]\label{defn:zderived}
An additive category 
$\tHMF{R}$ of graded matrix factorizations is defined 
by the following data. 
The set of objects is given by the set of 
all graded matrix factorizations: 
\begin{equation*}
\Ob(\tHMF{R}):= \tMF{R}\ .
\end{equation*}
For any two objects $\ti{M},\ti{M}'\in\tMF{R}$, 
the set of morphisms is given by
\begin{equation*}
\HomDz{R}(\ti{M},\ti{M}'):=\Hom_\tMF{R}(\ti{M},\ti{M}')/\sim\ ,
\end{equation*}
where two elements $\Phi$, $\Phi'$ in $\Hom_\tMF{R}(\ti{M},\ti{M}')$ 
are equivalent $\Phi\sim\Phi'$ 
if there exists a {\em homotopy}, \ie, 
a pair $(h_0,h_1):(\ti{P}_0\to\ti{P}_1',\ti{P}_1\to\ti{P}_0')$ 
of graded $R$-homomorphisms 
such that $h_0$ is of degree minus two, $h_1$ is of degree zero and 
$\Phi'-\Phi=(\tau^{-h}(p_1')h_0+h_1p_0, p_0'h_1+\tau^h(h_0)p_1)$. 
The composition of morphisms is induced from that 
on $\Hom_\tMF{R}(\ti{M},\ti{M}')$. 
 \label{defn:gD}
\end{defn}
A graded matrix factorization is 
the zero object in $\tHMF{R}$ 
if and only if it is a direct sum of 
the graded matrix factorizations of the forms 
$(\mfs{\tau^n(R)}{1}{f}{\tau^n(R)})\in\tMF{R}$ and 
$(\mfs{\tau^{n'}(R)}{f\ }{1\ }{\tau^{n'+h}(R)})\in\tMF{R}$ 
for some $n,n'\in\Z$.
\begin{lem}
The category $\tHMF{R}$ is Krull-Schmidt, that is, \\
(a)\ for any two objects $\ti{M},\ti{M}'\in\tHMF{R}$, 
$\HomDz{R}(\ti{M},\ti{M}')$ is finite dimensional; \\
(b)\ for any object $\ti{M}\in\tHMF{R}$ and 
any idempotent $e\in\HomDz{R}(\ti{M},\ti{M})$, 
there exists a matrix factorization $\ti{M}'\in\tHMF{R}$ 
and a pair of morphisms 
$\Phi\in\HomDz{R}(\ti{M},\ti{M}')$, $\Phi'\in\HomDz{R}(\ti{M}',\ti{M})$ 
such that $e=\Phi'\Phi$ and $\Phi\Phi'=\Id_{\ti{M}'}$. 
 \label{lem:KS}
\end{lem}
\begin{pf}
(a)\ Due to Lemma \ref{lem:finite}, 
$\oplus_{n\in\Z}\HomDz{R}(\ti{M},\tau^n(\ti{M}'))$ 
is a finitely generated graded 
${R}/\left(\frac{\p f}{\p x},\frac{\p f}{\p y},
\frac{\p f}{\p z}\right)$-module. 
Since the Jacobi ring
${R}/\left(\frac{\p f}{\p x},\frac{\p f}{\p y},\frac{\p f}{\p z}\right)$ 
is finite dimensional, 
the space $\HomDz{R}(\ti{M},\ti{M}')$ 
is finite dimensional over $\C$. 

(b)\ Let $R_+$ be the maximal ideal of $R$ 
of all positive degree elements. 
Note that any graded matrix factorization is isomorphic 
in $\tHMF{R}$ to a graded matrix factorization 
whose entries belong to $\tau^n(R_+)$ for some $n\in\Z$. 
Thus, we may assume that 
$\ti{M}:=(\mfs{\ti{P}_0}{p_0}{p_1}{\ti{P}_1})\in\tHMF{R}$ 
is such a graded matrix factorization. 
Suppose that $\ti{M}$ has an idempotent 
$e\in\HomDz{R}(\ti{M},\ti{M})$, $e^2=e$. This implies that there exists 
$\hat{e}\in\Hom_{\tMF{R}}(\ti{M},\ti{M})$ such that 
\begin{equation}\label{KS}
 \hat{e}^2-\hat{e}=(\tau^{-h}(p_1)h_0+h_1p_0,p_0h_1+\tau^h(h_0)p_1)
\end{equation}
for some homotopy $(h_0,h_1)$ on $\HomDz{R}(\ti{M},\ti{M})$. 
However, since each entry of $p_0$ and $p_1$ belongs to 
$\tau^n(R_+)$, each entry in the right hand side also belongs to 
$\tau^n(R_+)$. Let 
$\pi:\Hom_{\tMF{R}}(\ti{M},\ti{M})\to\Hom_{\tMF{R}}(\ti{M},\ti{M})$ 
be the canonical projection given by restricting each entry on 
$R/R_+=\C$. Then, eq.(\ref{KS}) in fact implies that 
$\pi(\hat{e})^2-\pi(\hat{e})=0$. 
Thus, for $\pi(\hat{e})=:(\hat{e}_0,\hat{e}_1)$, 
defining a matrix factorization $\ti{M}'\in\tHMF{R}$ by 
\begin{equation*}
 \ti{M}':=\Big( 
 \mfl{\hat{e}_0\ti{P}_0}{\hat{e}_1p_0\hat{e}_0}
{\hat{e}_0p_1\hat{e}_1}{\hat{e}_1\ti{P}_1}
 \Big) \ ,
\end{equation*}
one obtains a pair of morphisms 
$\Phi\in\HomDz{R}(\ti{M},\ti{M}')$ and $\Phi'\in\HomDz{R}(\ti{M}',\ti{M})$ 
such that $e=\Phi'\Phi$ and $\Phi\Phi'=\Id_{\ti{M}'}$. 
\qed\end{pf}

One can see that $\tau$ induces an automorphism on $\tHMF{R}$, 
which we shall denote by the same notation
$\tau:\tHMF{R}\to\tHMF{R}$. 
Explicitly, the action of $\tau$ on 
$\ti{M}=(\mfs{\ti{P}_0}{p_0}{p_1}{\ti{P}_1})\in\tHMF{R}$ is given by 
\begin{equation*}
 \tau \Big(\mf{\ti{P}_0}{p_0}{p_1}{\ti{P}_1}\Big):=
 \Big(\mf{\tau(\ti{P}_0)}{\tau(p_0)}{\tau(p_1)}{\tau(\ti{P}_1)}\Big)\ .
\end{equation*}
The action of $\tau$ on morphisms are naturally induced from that 
on graded ${R}$-homomorphisms between two graded right ${R}$-modules. 

Also, we have the shift functor $T$ on $\HMF{R}$, 
the graded version of that in Definition \ref{defn:shift}.
\footnote{The shift functor $T$ is often denoted by $[1]$. }
\begin{defn}[Shift functor on $\tHMF{R}$]
The {\it shift functor} 
$T:\tHMF{R}\to\tHMF{R}$ is defined as follows. 
The action of $T$ on $\ti{M}\in\tHMF{R}$ is given by 
\begin{equation*}
 T\Big(\mf{\ti{P}_0}{p_0}{p_1}{\ti{P}_1}\Big)
 :=\Big(\mf{\ti{P}_1}{-p_1\quad}{-\tau^h(p_0)\quad}{\tau^h(\ti{P}_0)}\Big)\ .
\end{equation*}
For any $\ti{M},\ti{M}'\in\tHMF{R}$, the action of $T$ on 
$\Phi=(\phi_0,\phi_1)\in\HomDz{R}(\ti{M},\ti{M}')$ 
is given by 
\begin{equation*}
 T(\phi_0,\phi_1):=(\phi_1,\tau^h(\phi_0))\ . 
\end{equation*}
 \label{defn:shift_Z}
\end{defn}
We remark that the square $T^2$ of the shift functor is not isomorphic 
to the identity functor on $\tHMF{R}$ but $T^2=\tau^h$. 
\begin{defn}[Mapping cone]
For an element $\Phi=(\phi_0,\phi_1)\in\Hom_\tMF{R}(\ti{M},\ti{M}')$,  
the {\em mapping cone} $C(\Phi)\in\tMF{R}$ is defined by 
\begin{equation*}
 \begin{split}
 & C(\Phi):=\Big( \mf{C_0}{c_0}{c_1}{C_1}\Big)\ ,\text{where} \\
 & C_0:=\ti{P}_1\oplus\ti{P}_0'\ ,\qquad 
   C_1:=\tau^h(\ti{P}_0)\oplus\ti{P}_1'\ ,\qquad 
   c_0:=\bp -p_1 & 0 \\ \phi_1 & p_0'\ep\ ,\qquad 
   c_1:=\bp -\tau^h(p_0) & 0 \\ \tau^h(\phi_0) & p_1'\ep\ .
 \end{split}
\end{equation*}
 \label{defn:cone_Z}
\end{defn}
This mapping cone is well-defined. In fact, 
one can see that the degree of $c_0$ and $c_1$ are zero and two, 
since 
the graded ${R}$-homomorphisms $p_1:\ti{P}_1\to\ti{P}_0$ of degree two 
and $p_0:\ti{P}_0\to\ti{P}_1$ of degree zero 
induce graded ${R}$-homomorphisms $-p_1:\ti{P}_1\to\tau^h(\ti{P}_0)$ 
of degree zero and $-\tau^h(p_0):\tau^h(\ti{P}_0)\to\ti{P}_1$ 
of degree two, respectively. 

The following is stated in \cite{o:2} Proposition 3.4. 
\begin{thm}
The additive category $\tHMF{R}$ endowed 
with the shift functor $T$ and the distinguished triangles 
forms a triangulated category, 
where a distinguished triangle is defined by a sequence 
isomorphic to the sequence 
\begin{equation*}
 \ti{M}\stackrel{\Phi}{\to}\ti{M}' \to C(\Phi)\to T(\ti{M})\ 
\end{equation*}
for some $\ti{M},\ti{M}'\in\tMF{R}$ 
and $\Phi\in\Hom_\tMF{R}(\ti{M},\ti{M}')$. 
\end{thm}
\begin{pf}
As in the case for $\HMF{A}$, 
the proof is the same as the proof of the analogous 
result for a usual homotopic category 
(see e.g. \cite{gm:1}, \cite{ks:1}). 
\qed\end{pf}
Let $\ti{M}=(\mfs{\ti{P}_0}{p_0}{p_1}{\ti{P}_1})\in\tHMF{R}$ 
be a graded matrix factorization of size $r$. 
One can choose homogeneous free basis 
$(b_1,\cdots,b_r;\bb_1,\cdots,\bb_r)$ such that 
$\ti{P}_0=b_1{R}\oplus\cdots\oplus b_r{R}$ and 
$\ti{P}_1=\bb_1{R}\oplus\cdots\oplus\bb_r{R}$. 
Then, the graded matrix factorization $\ti{M}$ is expressed 
as a pair $(Q,S)$ of $2r$ by $2r$ matrices, 
where $S$ is the diagonal matrix of the form 
$S:=\diag(s_1,\cdots,s_r; \sb_1,\cdots,\sb_r)$ such that 
$s_i=\deg(b_i)$ and $\sb_i=\deg(\bb_i)-1$ 
for $i=1,2,\cdots,r$ and 
\begin{equation}\label{Q}
 Q= \bp \0 & \varphi \\ \psi & \0 \ep\ ,\qquad \varphi,\psi\in\Mat_r(R)
\end{equation}
satisfying 
\begin{equation}
 Q^2=f\cdot\1_{2r}\ ,\qquad -SQ+QS+2EQ= Q\ .
 \label{S}
\end{equation}
%AAA
We call this $S$ a {\em grading matrix} of $Q$. 
This procedure $\ti{M}\mapsto (Q,S)$ 
gives the equivalence between 
the triangulated category $\tHMF{R}$ and 
the triangulated category $D^b_\Z(\A_f)$ in \cite{t:2}. 
This implies that 
$\tHMF{R}$ is an enhanced triangulated category 
in the sense of Bondal-Kapranov \cite{bk:1}. 

We shall represent the matrix factorization 
$\ti{M}=(\mfs{\ti{P}_0}{p_0}{p_1}{\ti{P}_1})$ 
by $(Q,S)$.

%%%%%%%%%%%%%%%%%%%%%%%%%%%%%%%%%%%%%%%%%%%%%%%%%%%%%%%%%%%%%%%%%%%%%%%%%%%%%%
\section{$\tHMF{R}$ for type ADE and representations of Dynkin quivers}
\label{sec:main}

In this section, we formulate the main theorem (Theorem \ref{thm:main}) 
of the present paper 
in subsection \ref{ssec:main-state}. 
The proof of the theorem is given 
in subsection \ref{ssec:proof}.

\subsection{Statement of the main theorem (Theorem \ref{thm:main})}\hfill\\
\label{ssec:main-state}
The main theorem states an equivalence 
between the triangulated category 
$\tHMF{R}$ for a polynomial $f\in R$ of type ADE 
with the derived category of modules over 
a {\em path algebra} of a {\em Dynkin quiver}. 
In order to formulate the results, 
we recall 
(i) the weighted homogeneous polynomials of type ADE and 
(ii) the notion of the path algebras of the Dynkin quivers. 

\vspace*{0.2cm}

\noindent
(i)\ {\bf ADE polynomials.}\ \, 
For a regular weight system $W$, we have the following facts 
\cite{sa:weight}. 

\vspace*{0.15cm}

\noindent
(a)\ 
There exist integers $m_1,\cdots,m_l$, 
called the {\em exponents} of $W$, such that 
$\chi_W(T)=T^{m_1}+\cdots + T^{m_l}$, where the smallest 
exponent is given by $\epsilon:=a+b+c-h$.

\vspace*{0.15cm}

\noindent
(b)\ 
The regular weight systems with $\epsilon>0$ are 
listed as follows. 
\begin{equation}
 \begin{split}
 & A_l : (1,b,l+1-b; l+1)\ ,1\le b\le l\ , \quad 
 D_l : (l-2,2,l-1;2(l-1))\ , \\ 
 & E_6 : (4,3,6;12)\ , \quad 
 E_7 : (6,4,9;18)\ , \quad
 E_8 : (10,6,15;30) \ .
 \end{split}
 \label{ADE}
\end{equation}
Here, the naming in the left hand side is given 
according to the identifications of the exponents of 
the weight systems with those of the simple Lie algebras.  
As a consequence, 
one obtains $\epsilon=1$ for all regular weight 
system with $\epsilon>0$. 
For the polynomials of type ADE, 
without loss of generality we may choose the followings: 
\begin{equation*}
f(x,y,z)=
\left\{\begin{array}{lll}
x^{l+1}+yz,      & h=l+1,    & A_{l}\ ( l\ge 1),\\
x^2y+y^{l-1}+z^2,& h=2(l-1), & D_l\  ( l\ge 4),\\
x^3+y^4+z^2,     & h=12,     & E_6,\\
x^3+xy^3+z^2,    & h=18,     & E_7, \\
x^3+y^5+z^2,     & h=30,     & E_8\ .
\end{array}\right.
\end{equation*}

\vspace*{0.25cm}

\noindent
\ (ii)\ {\bf Path algebras.}\ \, 
(a)\ The {\em path algebra} $\C\vec{\Delta}$ of a {\em quiver} 
is defined as follows
(see \cite{g:1}, \cite{r:1} and \cite{h} Chapter 1, 5.1). 
A {\em quiver} $\vec{\Delta}$ is a pair $(\Delta_0,\Delta_1)$ of 
the set $\Delta_0$ of vertices and 
the set $\Delta_1$ of arrows (oriented edges). 
Any arrow in $\Delta_1$ has a starting point and end point in $\Delta_0$. 
A path of {\em length} $r\ge 1$ 
from a vertex $v$ to a vertex $v'$ in a quiver $\vec{\Delta}$ 
is of the form $(v|\alpha_1,\cdots,\alpha_r|v')$ with 
arrows $\alpha_i\in\Delta_1$ satisfying 
the starting point of $\alpha_1$ is $v$, 
the end point of $\alpha_i$ is equal to the starting point of 
$\alpha_{i+1}$ 
for all $1\le i\le r-1$, 
and the end point of $\alpha_r$ is $v'$. 
In addition, we also define 
a path of length zero $(v|v)$ 
for any vertex $v$ in $\vec{\Delta}$. 
The {\em path algebra} $\C\vec{\Delta}$ of a quiver $\vec{\Delta}$ 
is then the $\C$-vector space with basis the set of all paths 
in $\vec{\Delta}$. 
The product structure is defined by the composition of paths, 
where the product of two non-composable paths is set to be zero. 

The category of 
finitely generated right modules over the path algebra 
$\C\vec{\Delta}$ is denoted by $\mod\C\vec{\Delta}$. 
It is an abelian category, and its derived category is 
denoted by $D^b(\mod\C\vec{\Delta})$. 
If $\Delta_0$ and $\Delta_1$ are finite sets 
and $\vec{\Delta}$ does not have any 
oriented cycle, then $\C\vec{\Delta}$ is a finite dimensional 
algebra and $D^b(\mod\C\vec{\Delta})$ 
is a Krull-Schmidt category.

\vspace*{0.15cm}

\noindent
(b)\ 
A {\em Dynkin quiver} $\vec{\Delta}$ of type ADE 
is one of the Dynkin diagram $\Delta$ 
listed in Figure 1 %\ref{fig:Dynkin} 
together with an orientation 
for each edge of the diagram. 
It is known \cite{g:1} that the number of all the isomorphism classes of
the indecomposable objects of 
the abelian category $\mod\C\vec{\Delta}$ of a quiver 
$\vec{\Delta}$ is finite if and only if the quiver $\vec{\Delta}$ is
a Dynkin quiver (of type ADE). 
\begin{equation*}
 \begin{split}
  A_l :\ &\qquad
\xymatrix{
 \bullet_1 \ar@{-}@{-}[r]^{}
  & \bullet_2 \ar@{-}[r] 
   & \, \cdots\, \ar@{-}[r]^{} 
   & \bullet_b \ar@{-}[r]^{}
   & \, \cdots\cdot\cdot\, \ar@{-}[r]^{}
  & \bullet_{l-1}  \ar@{-}[r]^{} 
  & \bullet_l \quad ,
}\\
%\end{equation*}
%\begin{equation*}
% \begin{array}{c} \\  
 D_l:\ 
%\end{array}
&\qquad
\xymatrix{
 & & & & \bullet_{l-1} \\
 \bullet_1 & 
 \bullet_2 \ar@{-}[l]^{} \ar@{-}[r]^{} &
 \cdot\cdots\cdot \ar@{-}[r] & 
 \bullet_{l-3} \ar@{-}[r] & 
 \bullet_{l-2} \ar@{-}[r]^{} \ar@{-}@<1.3ex>[u]^{} & 
 \bullet_l \quad ,
}\\
%\end{equation*}
%\begin{equation*}
 E_6:\ &\qquad
\xymatrix{
 & & \bullet_1 \ar@{-}@<-0.4ex>[d]^{} & \\
\bullet_5 & \bullet_3 \ar@{-}[l]^{} & 
\bullet_2 \ar@{-}[l] \ar@{-}[r] 
 & \bullet_4 \ar@{-}[r]^{} & \bullet_6\quad , 
} \\
%\end{equation*}
%\begin{equation*}
 E_7 :\ &\qquad 
\xymatrix{
 & & \bullet_5 & & \\
 \bullet_7 
 & \bullet_6 \ar@{-}[l]
 & \bullet_4 \ar@{-}[l]^{}\ar@{-}@<0.4ex>[u]^{} 
 & \bullet_3 \ar@{-}[l]
 & \bullet_2 \ar@{-}[l]^{} 
 & \bullet_1 \ar@{-}[l]\quad ,
} \\
%\end{equation*}
%\begin{equation*}
 E_8 :\ &\qquad  
\xymatrix{
 & & & & \bullet_8 \\
 \bullet_1 \ar@{-}[r]^{} 
 & \bullet_2 \ar@{-}[r] 
 & \bullet_3 \ar@{-}[r]^{} 
 & \bullet_4 \ar@{-}[r] 
 & \bullet_5 \ar@{-}[r]^{}\ar@{-}@<0.4ex>[u]^{}
 & \bullet_6 \ar@{-}[r] 
 & \bullet_7 \quad .
} 
 \end{split}
\end{equation*}
\begin{center}
 F{\scriptsize{IGURE}} 1.\ ADE Dynkin diagram
\end{center}
%\begin{figure}[h]
% \caption{ADE Dynkin diagram}
% \label{fig:Dynkin}
%\end{figure}
For a Dynkin diagram $\Delta$, 
we shall denote by $\Pi$ the set of vertices. 
For later convenience, 
the elements of $\Pi$ are labeled by the integers $\{1,\cdots, l\}$ 
as in the above figures, 
and we shall sometimes confuse vertices in $\Pi$ with 
labels in $\{1,\cdots, l\}$.

The following is the main theorem of the present paper. 
\begin{thm}
Let $f\in\C[x,y,z]$ be a polynomial of type ADE, 
and let $\vec{\Delta}$ be a Dynkin quiver 
of the corresponding type with a fixed orientation. 
Then, we have the following equivalence of the triangulated categories
\begin{equation*}
 \tHMF{R}\simeq D^b(\mod\C\vec{\Delta})\ .
\end{equation*}
 \label{thm:main}
\end{thm}

%%%%%%%%%%%%%%%%%%%%%%%%%%%%%%%%%%%%%%%%%%%%%%%%%%%%%%%%%%%%%%%%%%%%%%%%%%%%
\subsection{The proof of Theorem \ref{thm:main}}\hfill\\
 \label{ssec:proof}
The construction of the proof of Theorem \ref{thm:main} is 
as follows. 
\begin{itemize}
 \item[Step 1.]\ 
We describe the Auslander-Reiten (AR-)quiver 
for the triangulated category $\HMF{\Oh}$ 
of matrix factorizations 
due to \cite{e:1}, \cite{ar:1} and \cite{a:1} and 
give the matrix factorizations explicitly.  

 \item[Step 2.]\ We determine the structure of 
the triangulated category $\tHMF{R}$ of graded matrix factorizations 
(Theorem \ref{thm:grading}).

 \item[Step 3.]\ 
By comparing the AR-quiver of the 
category $\HMF{\Oh}$ with the category $\tHMF{R}$ 
we find the exceptional collections corresponding to
$\vec{\Delta}$ in $\tHMF{R}$ 
and complete the proof of the main theorem (Theorem \ref{thm:main}). 
\end{itemize}

\noindent
{\bf Step 1.\ \ The Auslander-Reiten quiver for $\HMF{\Oh}$. }

We recall the known results on the 
equivalence of the McKay quiver for Kleinean group and 
the AR-quiver for the simple singularities \cite{ar:1}, 
\cite{a:1} and \cite{e:1}. 

Let $\CC$ be a Krull-Schmidt category over $\C$ 
which is equivalent to a small category. 
%For a Krull-Schmidt category $\CC$ 
%over $\C$, a
An object $X\in\Ob(\CC)$ is called {\em indecomposable} 
if any idempotent 
$e\in\Hom_\CC(X,X)$ is zero or the identity $\Id_X$. 
For two objects 
$X,Y\in\Ob(\CC)$, denote by 
$\rad_\CC(X,Y)$ the linear subspace of 
$\Hom_\CC(X,Y)$ of non-invertible morphisms 
from $X$ to $Y$. 
We denote by $\rad_\CC^2(X,Y)\subset\rad_\CC(X,Y)$ 
the space of morphisms each of which is described 
as a composition $\Phi'\Phi$ with $\Phi\in\rad_\CC(X,Z)$, 
$\Phi'\in\rad_\CC(Z,Y)$ for some object $Z\in\Ob(\CC)$. 
For two indecomposable objects $X,Y\in\Ob(\CC)$, 
an element in $\rad_\CC(X,Y)\backslash\rad_\CC^2(X,Y)$ 
is called an {\em irreducible morphism}. 
The space 
$\Irr_\CC(X,Y):=\rad_\CC(X,Y)/\rad_\CC^2(X,Y)$ in fact forms 
a subvector space of $\Hom_\CC(X,Y)$. 
We call 
by the {\em AR-quiver $\Gamma(\CC)$} of a Krull-Schmidt category $\CC$ 
the quiver $\Gamma(\CC):=(\Gamma_0,\Gamma_1)$ 
whose vertex set $\Gamma_0$ consists of 
the isomorphism classes $[X]$ of the indecomposable objects 
$X\in\Ob(\CC)$ and 
whose arrow set $\Gamma_1$ consists of 
$\dim_\C(\Irr_\CC(X,Y))$ arrows 
from $[X]\in\Gamma_0$ to $[Y]\in\Gamma_0$ 
for any $[X],[Y]\in\Gamma_0$ (see \cite{r:1},\cite{h},\cite{y}).

On the other hand, 
for a Dynkin diagram $\Delta$ listed in Figure 1, 
we define a quiver consisting of 
the vertex set $\Pi$ and 
arrows in both directions 
$\mfss{k}{}{}{k}\!'$ for each edge 
of $\Delta$ between vertices $k,k'\in\Pi$. 
The resulting quiver is denoted by $\lrDelta$. 

Note that the category $\HMF{\Oh}$ 
is Krull-Schmidt (see \cite{y}, Proposition 1.18). 
\begin{thm}
Let $f$ be a polynomial of type ADE, which we regard 
as an element of $\Oh$. 

{\em (i)}\ $($\cite{ar:1},\cite{a:1},\cite{e:1}$)$
The AR-quiver of the category $\HMF{\Oh}$ 
is isomorphic to the quiver $\lrDelta$ 
corresponding to $f$ $:$
\begin{equation*}
 \Gamma(\HMF{\Oh})\simeq \lrDelta\ .
\end{equation*}

{\em (ii)}\ 
According to (i), 
fix an identification $\Pi\simeq\Gamma_0(\HMF{\Oh})$, 
$k\mapsto [M^k]$. 
A representative $M^k$ of the isomorphism classes $[M^k]$ 
of the indecomposable matrix factorizations 
of minimum size is given explicitly in Table \ref{tbl:1}. 
The size of $M^k$ is $2\nu_k$, 
where $\nu_k$ is the coefficient 
of the highest root for $k\in\Pi$. 
 \label{thm:AR}
\end{thm}
\begin{pf}
(i)\ This statement follows from the combination of results of
\cite{ar:1} and \cite{e:1}, 
where the Auslander-Reiten quivers of the categories of 
the maximal Cohen-Macaulay modules over $\Oh/(f)$ 
for type ADE are determined in \cite{ar:1}, 
and the equivalence of the category of Maximal Cohen-Macaulay modules 
with the category $\HMF{\Oh}$ of the matrix factorizations is 
given in \cite{e:1}. 

(ii)\ 
Since we have $\sharp(\Pi)$ non-isomorphic 
matrix factorizations (Table \ref{tbl:1}), 
these actually complete all the vertices of the AR-quiver 
$\Gamma(\HMF{\Oh})$. 
\qed \end{pf}
\begin{rem}
In \cite{y}, matrix factorizations for a polynomial of type ADE 
in two variables $x,y$ are listed up completely. 
On the other hand, 
for type $A_l$ and $D_l$ in both two and three variables, 
all the matrix factorizations and the AR-quivers 
are presented in \cite{sch}, 
where the relation of the results in two variables and those in three
variables is given. This gives a method of finding 
the matrix factorizations of a polynomial of type $E_l$, $l=6,7,8$, 
in three variables 
from the ones in two variables case \cite{y}. 
For a recent paper in physics, see also \cite{kl:3}. 
 \label{rem:mf}
\end{rem}
Hereafter we fix an identification of 
$\Gamma(\HMF{\Oh})$ with $\lrDelta$ 
by $k\leftrightarrow [M^k]$.

\vspace*{0.3cm}

\noindent
{\bf Step 2.\ \ Indecomposable objects in $\tHMF{R}$}

Recall that one has the inclusions $R\subset\O\subset\Oh$. 
We prepare some definitions 
for any fixed weighted homogeneous polynomial $f\in R$. 
\begin{defn}[Forgetful functor from $\tHMF{R}$ to $\HMF{\Oh}$]
For a fixed weighted homogeneous polynomial $f\in R$, 
there exists a functor $F:\tHMF{R}\to\HMF{\Oh}$ 
given by $F(\ti{M}):=\ti{M}\otimes_R\Oh$ 
for $\ti{M}\in\tHMF{R}$ and 
the naturally induced homomorphism 
$F:\HomDz{R}(\ti{M},\ti{M}')\to\HomD{\Oh}(F(\ti{M}),F(\ti{M}'))$ 
for any two objects $\ti{M}, \ti{M}'\in\tHMF{R}$. 
We call this $F$ the {\em forgetful functor}. 
 \label{defn:F}
\end{defn}
It is known that $F:\tHMF{R}\to\HMF{\Oh}$ 
brings an indecomposable object to 
an indecomposable object (\cite{y}, Lemma 15.2.1).

Let us introduce the notion of {\em distance} between 
two indecomposable objects in $\HMF{\Oh}$ and in $\tHMF{R}$ as follows. 
\begin{defn}[Distance]
For any two indecomposable objects $M,M'\in\HMF{\Oh}$, 
define the {\em distance} $d(M,M')\in\Z_{\ge 0}$ from $M$ to $M'$ by 
the minimal length of the paths from $[M]$ to $[M']$ 
in the AR-quiver $\Gamma(\HMF{\Oh})$ of $\HMF{\Oh}$. 
In particular, we have $d(M,M')=0$ if and only if $M\simeq M'$ 
in $\HMF{\Oh}$.

For any two indecomposable objects $\ti{M},\ti{M}'\in\tHMF{R}$, 
the {\em distance} $d(\ti{M},\ti{M}')\in\Z_{\ge 0}$ 
from $\ti{M}$ to $\ti{M}'$ is defined by 
\begin{equation*}
 d(\ti{M},\ti{M}'):=d(F(\ti{M}),F(\ti{M}'))\ .
\end{equation*}
 \label{defn:distance}
\end{defn}
By definition, 
an irreducible morphism 
exists in $\HomD{\Oh}(M,M')$ if and only if $d(M,M')=1$. 

Let us return to the case that $f$ is of type ADE. 
In this case,
the distance is in fact symmetric: 
$d(M^k,M^{k'})=d(M^{k'},M^k)$ for any $[M^k],[M^{k'}]\in\Pi$, 
since, due to Theorem \ref{thm:AR}, there exists an arrow from 
$[M^k]$ to $[M^{k'}]$ if and only if 
there exists an arrow from $[M^{k'}]$ to $[M^k]$. 
For two indecomposable objects $\ti{M}^k,\ti{M}^{k'}\in\tHMF{R}$ 
such that $F(\ti{M}^k)=M^k$, $F(\ti{M}^{k'})=M^{k'}$, 
we denote $d(\ti{M}^k,\ti{M}^{k'})=d(M^k,M^{k'}):=d(k,k')$.

In order to state the following Theorem \ref{thm:grading}, 
it is convenient 
to introduce the ordered decomposition $\Pi=\{\Pi_1,\Pi_2\}$ 
of $\Pi$, called 
a {\em principal decomposition} of $\Pi$ \cite{sa:principal}, 
as follows. 
We first define the base vertex $[M_o]\in\Pi$. 
For a diagram of $D_l$ or $E_l$, $l=6,7,8$, 
we choose $[M_o]$ as the trivalent vertex. % on which three edges join. 
Explicitly, it is $[M^{k=l-2}]$ for type $D_l$, 
$[M^{k=2}]$ for type $E_6$, 
$[M^{k=4}]$ for type $E_7$ and $[M^{k=5}]$ for type $E_8$ (Figure 1). 
For type $A_l$, we set $[M_o]:=[M^{k=b}]$ 
%(see Figure \ref{fig:Dynkin}). 
(see Figure 1) depending on the index $b$, $1\le b\le l$ 
(see eq.(\ref{ADE})). 
Then, we define 
the decomposition of $\Pi$ by 
\begin{equation*}
 \Pi_1:=\{k\in\Pi\ |\ d(M_o,M^k)\in 2\Z_{\ge 0}+1\}\ ,\qquad 
 \Pi_2:=\{k\in\Pi\ |\ d(M_o,M^k)\in 2\Z_{\ge 0}\}\ .
\end{equation*}

Recall that we express a graded matrix factorization 
$\ti{M}\in\tHMF{R}$ by the pair $\ti{M}=(Q,S)$, 
where $Q$ denotes a matrix factorization in eq.(\ref{Q}) and 
$S$ is the grading matrix defined in eq.(\ref{S}). 
\begin{thm} Let $f\in R$ be a polynomial of type ADE. 
The triangulated category $\tHMF{R}$ 
is described as follows. 

\noindent
{\em (i)}\ {\bf (Objects)}$:$\ \, 
The set of isomorphism classes of all indecomposable objects of 
$\tHMF{R}$ is given by 
\begin{equation*}
[\ti{M}^k_n:=(Q^k,S^k_n)],
\quad k\in\Pi\ ,
\quad n\in\Z \ .
\end{equation*}
Here, $\bullet$\ 
$Q^k$ is the matrix factorizations of size $2\nu_k$ 
given in Theorem \ref{thm:AR} (ii) (Table \ref{tbl:1}), 

$\bullet$\ 
the grading matrix $S^k_n$ for $k\in\Pi_\sigma$, $\sigma=1,2$,  
and $n\in\Z$ is given by$:$ 
\begin{equation}
 S^k_n:
=\diag\left(q^k_1,-q^k_1,\cdots,q^k_{\nu_k},-q^k_{\nu_k};
 \qb^k_1,-\qb^k_1,\cdots,\qb^k_{\nu_k},-\qb^k_{\nu_k}\right)
+\phi^k_n\cdot\1_{4\nu_k}\ ,
 \label{gradADE}
\end{equation}
where the data of the first term, called the traceless part$:$ 
$q^k_j\in\frac{2}{h}\Z-\frac{\sigma}{h}$ and 
$\qb^k_j\in\frac{2}{h}\Z-\frac{\sigma}{h}-1$ for $1\le j\le\nu_k$ 
are given in {\em Table \ref{tbl:2}}, and 
the coefficient of the second term, called the phase, 
is given by $\phi^k_n:=\phi(\ti{M}^k_n)=\frac{2n+\sigma}{h}$.

\vspace*{0.1cm}

\noindent
{\em (ii)}\ {\bf (Morphisms)}$:$\ \, 
{\em (ii-a)}\ An irreducible morphism exists in 
$\HomDz{R}(\ti{M}^k_n,\ti{M}^{k'}_{n'})$ if and only if 
$d(k,k')=1$ and 
$\phi^{k'}_{n'}=\phi^k_n+\ov{h}$.

{\em (ii-b)}\ 
$\dim_\C(\HomDz{R}(\ti{M}^k_n,\ti{M}^{k'}_{n'}))=1$ for 
$\phi^{k'}_{n'}-\phi^k_n=\ov{h} d(k,k')$.

\vspace*{0.1cm}

\noindent
{\em (iii)}\ {\bf (The Serre duality)}$:$\ \, 
The automorphism $\S:=T\tau^{-1}:\tHMF{R}\to\tHMF{R}$ 
satisfies the following properties. 

{\em (iii-a)}\ 
$\HomDz{R}(\ti{M}^k_n,\S(\ti{M}^k_n))\simeq\C$ 
for any indecomposable object $\ti{M}^k_n\in\tHMF{R}$. 

{\em (iii-b)}\ 
This isomorphism of (iii-a) induces the following 
bilinear map $:$ 
\begin{equation*}
 \HomDz{R}(\ti{M}^k_n,\ti{M}^{k'}_{n'})
 \otimes\HomDz{R}(\ti{M}^{k'}_{n'},\S(\ti{M}^k_n))\to\C\ ,
\end{equation*}
which is nondegenerate for any $k,k'\in\Pi$ and $n,n'\in\Z$. 
 \label{thm:grading}
\end{thm}
\begin{pf}
(i)\ 
For each $M^k$, 
by direct calculations based on the explicit form of $M^k$ in 
Table \ref{tbl:1}, 
we can attach grading matrices $S$ satisfying eq.(\ref{S}). 
This, together with the fact that 
$F:\tHMF{R}\to\HMF{\Oh}$ brings an indecomposable object 
to an indecomposable object, implies that 
the union $\coprod_{k=1}^lF^{-1}(M^k)$ gives the set of all 
indecomposable objects in $\tHMF{R}$ and then 
$F:\tHMF{R}\to\HMF{\Oh}$ is a surjection 
(see also \cite{ar:2} and \cite{y}, Theorem 15.14). 
Actually, $S$ is unique up to an addition of 
a constant multiple of the identity. 
Therefore, we decompose $S$ into the traceless part and the phase
part. 
%as in eq.(\ref{gradADE}). 
Due to the restriction of degrees (\ref{P-grad}) 
and the definition of $S$, 
one has $\pm q_j^k+\phi^k_n\in\frac{2}{h}\Z$ and 
$\pm \qb_j^k+\phi^k_n\in\frac{2}{h}\Z-1$ for any 
$k\in\Pi$ and $1\le j\le \nu_k$. 
% (see below eq.(\ref{S})), 
%and b
By solving these conditions, we obtain 
Statement (i), \ie, Table \ref{tbl:2}. 
In particular, one has 
$F^{-1}(M^k)=\{\ti{M}^k_n\ |\ n\in\Z\}$ and 
$\tau(\ti{M}^k_n)=\ti{M}^k_{n+1}$.

\vspace*{0.3cm}

\noindent
(ii-a)\ 
For any two indecomposable objects 
$\ti{M}^k_n,\ti{M}^{k'}_{n'}\in\tHMF{R}$, 
Statement (i) implies that 
$F:\HomDz{R}(\ti{M}^k_n,\ti{M}^{k'}_{n'})\to\HomD{\Oh}(M^k,M^{k'})$ 
is injective and then 
$$F:\oplus_{n''\in\Z}\HomDz{R}(\ti{M}^k_n,\tau^{n''}(\ti{M}^{k'}_{n'}))
\simeq\HomD{\Oh}(M^k,M^{k'})$$ is an isomorphism of vector spaces. 
This induces the following isomorphisms: 
\begin{equation*}
 \begin{array}{cccc}
% F:&\oplus_{n''\in\Z}\HomDz{R}(\ti{M}^k_n,\tau^{n''}(\ti{M}^{k'}_{n'})) 
% &\overset{\sim}{\to} &\HomD{\Oh}(M^k,M^{k'})\ ,\\
%     &       \cup        &&  \cup \\
 F:&\oplus_{n''\in\Z}\rad_{\tHMF{R}}(\ti{M}^k_n,\tau^{n''}(\ti{M}^{k'}_{n'}))
 &\overset{\sim}{\to} &\rad_{\HMF{\Oh}}(M^k,M^{k'}) \ ,\\
     &       \cup        &&  \cup \\
 F:&
 \oplus_{n''\in\Z}\rad_{\tHMF{R}}^2(\ti{M}^k_n,\tau^{n''}(\ti{M}^{k'}_{n'}))
 &\overset{\sim}{\to}&\rad_{\HMF{\Oh}}^2(M^k,M^{k'}) \ ,
 \end{array}
\end{equation*}
and hence, the isomorphism 
$F:\oplus_{n''\in\Z}
\Irr_{\tHMF{R}}(\ti{M}^k_n,\tau^{n''}(\ti{M}^{k'}_{n'}))
\overset{\sim}{\to}\Irr_{\HMF{\Oh}}(M^k,M^{k'})$. 

For $k,k'\in\Pi$, 
define a multi-set $\fC(k,k')$ 
of non-negative integers by 
\begin{equation*}
\fC(k,k')
:=\{h(\phi(\tau^{n''}(\ti{M}^{k'}_{n'}))-\phi(\ti{M}^k_n))\ |\ 
 \HomDz{R}(\ti{M}^k_n,\tau^{n''}(\ti{M}^{k'}_{n'}))\ne 0\ ,\ \ n''\in\Z\}\ ,
\end{equation*} 
where we fix $n,n'\in\Z$ for each $F^{-1}(M^k)$ and the integer 
$h(\phi(\tau^{n''}(\ti{M}^{k'}_{n'}))-\phi(\ti{M}^k_n))
=h(\phi^{k'}_{n'}-\phi^k_n)+2n''$ 
appears with multiplicity 
$\dim_\C(\HomDz{R}(\ti{M}^k_n,\tau^{n''}(\ti{M}^{k'}_{n'})))$. 
The multi-set $\fC(k,k')$ in fact depends only on 
$k$ and $k'$, and is independent of the choice of $n$ and $n'$. 

For $k,k'\in\Pi$ such that $d(k,k')=1$, 
by calculating 
$\HomD{\Oh}(M^k,M^{k'})$ 
using the explicit forms of the matrix factorizations $M^k, M^{k'}$ 
in Table \ref{tbl:1}, 
one can see that $\fC(k,k')$ consists 
of positive odd integers including 
$1$ of multiplicity one. 
This implies that, 
for two indecomposable objects 
$\ti{M}^k_n,\ti{M}^{k'}_{n'}\in\tHMF{R}$, 
a morphism in $\HomDz{R}(\ti{M}^k_n,\ti{M}^{k'}_{n'})$ 
is an irreducible morphism if and only if 
$h(\phi^{k'}_{n'}-\phi^k_n)=1$ 
(Statement (ii-a)). 
(A morphism in $\HomDz{R}(\ti{M}^k_n,\ti{M}^{k'}_{n'})$ 
which corresponds to $c\in\fC(k,k')$ with $c\ge 3$ can be obtained 
by composing irreducible morphisms. )

\vspace*{0.3cm}

\noindent
(ii-b)\ 
One gets the following lemma: 
\begin{lem}\label{lem:A}
{\em 1)}\ For an indecomposable object $\ti{M}^k_n\in\tHMF{R}$, one has 
\begin{equation*}
 \text{\em 1-a)}\ \HomDz{R}(\ti{M}^k_n,\ti{M}^k_n)\simeq\C\ ,\qquad  
 \text{\em 1-b)}\ \HomDz{R}(\ti{M}^k_n,\S(\ti{M}^k_n))\simeq\C\ . 
\end{equation*}
{\em 2)}\ For an indecomposable object $\ti{M}^k_n\in\tHMF{R}$, 
a morphism $\Psi\in\HomDz{R}(\ti{M}^k_n,\S(\ti{M}^k_n))$ and 
an indecomposable object $\ti{M}^{k'}_{n'}\in\tHMF{R}$ such that 
$d(k,k')=1$, 

{\em 2-a)}\ $\Psi\Phi\sim 0$ holds for an irreducible morphism
$\Phi\in\HomDz{R}(\ti{M}^{k'}_{n'},\ti{M}^k_n))$, 

{\em 2-b)}\ $\S(\Phi)\Psi\sim 0$ holds for an irreducible morphism
$\Phi\in\HomDz{R}(\ti{M}^k_n,\ti{M}^{k'}_{n'})$. 
\end{lem}
\begin{pf}
By direct calculations of 
$\HomDz{R}(\ti{M}^k_n,\ti{M}^k_n)$ and 
$\HomDz{R}(\ti{M}^k_n,\S(\ti{M}^k_n))$ 
using the explicit forms of 
the matrix factorizations in Table \ref{tbl:1} again. 
\qed\end{pf}

Note that Lemma \ref{lem:A}.1-a also follows from the AR-quiver
$\Gamma(\HMF{\Oh})$ of $\HMF{\Oh}$ in Theorem \ref{thm:AR} 
and Statement (ii-a). 

Due to Lemma \ref{lem:A}.1-b, 
for any indecomposable object $\ti{M}^k_n\in\tHMF{R}$, 
one can consider the cone $C(\Psi)\in\tHMF{R}$ of 
a nonzero morphism $\Psi\in\HomDz{R}(\S^{-1}(\ti{M}^k_n),\ti{M}^k_n)$, 
that is, the cone of a lift of the 
morphism $\Psi\in\HomDz{R}(\S^{-1}(\ti{M}^k_n),\ti{M}^k_n)$ 
to a homomorphism in $\Hom_\tMF{R}(\S^{-1}(\ti{M}^k_n),\ti{M}^k_n)$. 
We denote it simply by $C(\Psi)$, since 
two different lifts of the morphism 
$\Psi\in\HomDz{R}(\S^{-1}(\ti{M}^k_n),\ti{M}^k_n)$ lead 
to two cones which are isomorphic to each other in $\tHMF{R}$. 

Recall that an {\em Auslander-Reiten (AR-)triangle} 
(see \cite{h},\cite{y}, and also an Auslander-Reiten sequence or 
equivalently an almost split sequence \cite{ar:ass}) 
in a Krull-Schmidt triangulated category 
$\CC$ is a distinguished triangle 
\begin{equation}\label{ARtri-general}
 X\overset{u}{\to} Y\overset{v}{\to} Z\overset{w}{\to} T(X)
\end{equation}
satisfying the following conditions: 
\begin{itemize}
 \item[(AR1)] $X, Z$ are indecomposable objects in $\Ob(\CC)$. 
 \item[(AR2)] $w\ne 0$
 \item[(AR3)] If $\Phi:W\to Z$ is not a split epimorphism, then 
there exists $\Phi':W\to Y$ such that $v\Phi'=\Phi$. 
\end{itemize}
Then, it is known (see \cite{h}, Proposition in 4.3) 
that $u$ and $v$ are irreducible morphisms. 
\begin{lem}\label{lem:ARtri}
For an indecomposable object $\ti{M}^k_n\in\tHMF{R}$ and 
the cone $C(\Psi)\in\tHMF{R}$ of 
a nonzero morphism $\Psi\in\HomDz{R}(\S^{-1}(\ti{M}^k_n),\ti{M}^k_n)$, 
the distinguished triangle 
\begin{equation}\label{ARtri}
 \ti{M}^k_n\to C(\Psi)\to\tau(\ti{M}^k_n)\overset{T(\Psi)}{\to} 
 T(\ti{M}^k_n)
\end{equation}
is an AR-triangle. 
\end{lem}
\begin{pf}
Since by definition the conditions (AR1) and (AR2) are already satisfied, 
it is enough to show that the condition (AR3) holds. 
Consider the functor $\HomDz{R}(\ti{M}^{k'}_{n'}, - )$ 
on the distinguished triangle (\ref{ARtri}) for 
an indecomposable object $\ti{M}^{k'}_{n'}\in\tHMF{R}$. 
One obtains the long exact sequence as follows: 
\begin{equation}\label{coh-functor}
 \begin{split}
 & \cdots\to\HomDz{R}(\ti{M}^{k'}_{n'},\S^{-1}(\ti{M}^k_n))
 \overset{\Psi}{\to}\HomDz{R}(\ti{M}^{k'}_{n'},\ti{M}^k_n) \\
 &\hspace*{3.0cm}\to\HomDz{R}(\ti{M}^{k'}_{n'},C(\Psi))\to \\
 &\qquad\qquad 
 \HomDz{R}(\ti{M}^{k'}_{n'},\tau(\ti{M}^k_n))\overset{T(\Psi)}{\to} 
 \HomDz{R}(\ti{M}^{k'}_{n'},T(\ti{M}^k_n))
 \to\cdots\ .
 \end{split}
\end{equation}
A nonzero morphism in $\HomDz{R}(\ti{M}^{k'}_{n'},\tau(\ti{M}^k_n))$ 
is a split epimorphism if and only if 
$\ti{M}^{k'}_{n'}=\tau(\ti{M}^k_n)$, in which case 
$\HomDz{R}(\ti{M}^{k'}_{n'},\tau(\ti{M}^k_n))\simeq\C$ is spanned 
by the identity 
$\Id_{\ti{M}^{k'}_{n'}}\in\HomDz{R}(\ti{M}^{k'}_{n'},\ti{M}^{k'}_{n'})$. 
If $\ti{M}^{k'}_{n'}\ne\tau(\ti{M}^k_n)$, 
due to Lemma \ref{lem:A}.2-a, the map $T(\Psi)$ in the long exact
sequence (\ref{coh-functor}) 
turns out to be zero, 
which implies that we have the surjection 
$\HomDz{R}(\ti{M}^{k'}_{n'},C(\Psi))
\to\HomDz{R}(\ti{M}^{k'}_{n'},\tau(\ti{M}^k_n))$ and 
then the condition (AR3). 
\qed\end{pf}
The existence of the AR-triangle together 
with the corresponding long exact sequence (\ref{coh-functor}) 
leads to the two key lemmas as follows. 
\begin{lem}
Let $\ti{M}^k_n\in\tHMF{R}$ be an indecomposable object. 
The cone $C(\Psi)$ of a nonzero morphism 
$\Psi\in\HomDz{R}(\S^{-1}(\ti{M}^k_n),\ti{M}^k_n)$ is 
isomorphic to the direct sum of indecomposable objects 
$\ti{M}^{k_1}_{n_1}\oplus\cdots\oplus\ti{M}^{k_m}_{n_m}$ 
for some $m\in\Z_{>0}$ 
such that $\{k_1,\cdots,k_m\}=\{k'\in\Pi\ |\ d(k,k')=1\}$ and 
$\phi(\ti{M}^{k_i}_{n_i})=\phi(\ti{M}^k_n)+\ov{h}$ 
for any $i=1,\cdots,m$. 
 \label{lem:B}
\end{lem}
\begin{pf}
For the AR-triangle in eq.(\ref{ARtri}), 
the morphisms 
$\ti{M}^k_n\to C(\Psi)$ and 
$C(\Psi)\to\tau(\ti{M}^k_n)$ are irreducible and hence 
one has $\phi(\ti{M}^{k_i}_{n_i})=\phi(\ti{M}^k_n)+\ov{h}$ 
and $d(k',k_i)=1$, 
for any $i=1,\cdots, m$, with the direct sum decomposition of 
indecomposable objects 
$C(\Psi)\simeq\oplus_{i=1}^m\ti{M}^{k_i}_{n_i}$ above. 
Also, the fact that the distinguished triangle (\ref{ARtri}) 
is an AR-triangle further guarantees that, 
for an indecomposable object $\ti{M}^{k'}_{n'}\in\tHMF{R}$, 
there exists an irreducible morphism in 
$\HomDz{R}(\ti{M}^{k}_{n},\ti{M}^{k'}_{n'})$ 
if and only if $(k',n')=(k_i,n_i)$ for some $i=1,\cdots, m$ 
(see Lemma in p.40 of \cite{h}). 
Thus, the rest of the proof is to show $k_i\ne k_j$ if $i\ne j$, 
for which it is enough to check this lemma 
at the level of the grading matrices 
(see eq.(\ref{cone-grading}) below Table \ref{tbl:2}). 
\qed\end{pf}
\begin{lem}
For an indecomposable object $\ti{M}^{k}_{n}\in\tHMF{R}$, 
let $C(\Psi)$ be the cone of a nonzero morphism 
$\Psi\in\HomDz{R}(\S^{-1}(\ti{M}^{k}_{n}),\ti{M}^{k}_{n})$, 
and let 
$C(\Psi)\simeq\oplus_{i=1}^m\ti{M}^{k_i}_{n_i}$, 
$m\in\Z_{>0}$, 
be the direct sum decomposition of indecomposable objects as above. 
Then, for an indecomposable object 
$\ti{M}^{k'}_{n'}\in\tHMF{R}$ and 
the given multi-set $\fC(k',k)$, one obtains 
\begin{equation*}
 \coprod_{i=1}^m\fC(k',k_i)
 =\{c-1\ |\ c\in\fC(k',k), c\ne 0\}
 \coprod 
 \{c+1\ |\ c\in\fC(k',k), c\ne h-2\ \text{if}\ k={k'}^\S\ \}\ ,
\end{equation*}
where, for $k\in\Pi$, 
$k^\S\in\Pi$ denotes the vertex such that 
$[F(\S(\ti{M}^k_n))]=[M^{k^\S}]\in\Pi$ for any $n\in\Z$. 
 \label{lem:C}
\end{lem}
\begin{pf}
Consider the long exact sequence (\ref{coh-functor}) for the
AR-triangle (\ref{ARtri}). 
As discussed in the proof of Lemma \ref{lem:ARtri}, 
if $\ti{M}^{k'}_{n'}\ne\tau(\ti{M}^k_n)$, 
the map $T(\Psi)$ in eq.(\ref{coh-functor}) 
is zero due to Lemma \ref{lem:A}.2-a. 
Similarly, if $\ti{M}^{k'}_{n'}\ne\S^{-1}(\ti{M}^k_n)$, 
the map $\Psi$ in eq.(\ref{coh-functor}) 
is zero due to Lemma \ref{lem:A}.2-b. 
On the other hand, if $\ti{M}^{k'}_{n'}=\tau(\ti{M}^k_n)$, 
the map $T(\Psi)$ in eq.(\ref{coh-functor}) is 
an isomorphism and hence the map 
$\HomDz{R}(\ti{M}^{k'}_{n'},C(\Psi))
\to\HomDz{R}(\ti{M}^{k'}_{n'},\tau(\ti{M}^k_n))$ turns out to be zero. 
Similarly, if $\ti{M}^{k'}_{n'}=\S^{-1}(\ti{M}^k_n)$, 
the map $\Psi$ in eq.(\ref{coh-functor}) is 
an isomorphism and hence the map 
$\HomDz{R}(\ti{M}^{k'}_{n'},\ti{M}^k_n)
\to\HomDz{R}(\ti{M}^{k'}_{n'},C(\Psi))$ 
turns out to be zero. 
Combining these facts, 
one has the exact sequences as follows: 
\begin{itemize}
 \item $0\to\Hom(\ti{M}^k_n,\ti{M}^{k'}_{n'})
\to\HomDz{R}(\ti{M}^k_n,C(\Psi))
\to\Hom(\ti{M}^k_n,\tau(\ti{M}^{k'}_{n'}))\to 0$ 
if $\ti{M}^k_n\ne\tau(\ti{M}^{k'}_{n'})$ and 
$\ti{M}^k_n\ne\S^{-1}(\ti{M}^{k'}_{n'})$,

 \item 
 $0\to\Hom(\ti{M}^k_n,\ti{M}^{k'}_{n'})
 \to\HomDz{R}(\ti{M}^k_n,C(\Psi))\to 0$ 
if $\ti{M}^k_n=\tau(\ti{M}^{k'}_{n'})$, and 

 \item 
 $0\to\HomDz{R}(\ti{M}^k_n,C(\Psi))
\to\Hom(\ti{M}^k_n,\tau(\ti{M}^{k'}_{n'}))\to 0$ 
if $\ti{M}^k_n=\S^{-1}(\ti{M}^{k'}_{n'})$, though in this case 
$\Hom(\ti{M}^k_n,\tau(\ti{M}^{k'}_{n'}))=0$. 
\end{itemize}
{}From these results follows the statement of this lemma. 
\qed\end{pf}
Calculate $\HomDz{R}(M^k,\tau^n(M^k))$, $n\in\Z$, 
where we set $k=1$ for type $A_l$, $D_l$ and $E_8$, $k=5$ or $6$ for $E_6$, 
and $k=7$ for $E_7$, 
by using the explicit forms of matrix factorizations $M^k$ 
in Table \ref{tbl:1}. 
Then, one obtains $\fC(k,k)$. 
Using Lemma \ref{lem:B} and Lemma \ref{lem:C} repeatedly, 
one can actually obtain $\fC(k',k'')$ for all $k',k''\in\Pi$ 
(Table \ref{tbl:3}). 
In particular, by a use of Table \ref{tbl:3}, 
one can get Statement (ii-b) and 
\begin{cor}
The following equivalent statements hold: for two indecomposable
objects $\ti{M}^k_n,\ti{M}^{k'}_{n'}\in\tHMF{R}$, 

{\em (a)}\ If 
$\phi(\S(\ti{M}^{k'}_{n'}))-\phi(\ti{M}^{k}_{n})< 
\ov{h}d(k,{k'}^\S)$, 
then $\HomDz{R}(\ti{M}^{k'}_{n'},\ti{M}^{k}_{n})=0$.

{\em (b)}\ If 
$\phi(\ti{M}^{k}_{n})-\phi(\ti{M}^{k'}_{n'})
<\ov{h}d(k',k)$, 
then $\HomDz{R}(\ti{M}^{k}_{n},\S(\ti{M}^{k'}_{n'}))=0$. 
 \label{cor:serre1}
\qed\end{cor}

\vspace*{0.3cm}

\noindent
(iii-a)\ 
This is already proven in Lemma \ref{lem:A}.1-b. 

\vspace*{0.3cm}

\noindent
(iii-b)\ 
Suppose that there exists a nonzero morphism 
$\Phi\in\HomDz{R}(\ti{M}^{k'}_{n'},\ti{M}^{k}_{n})$ 
for two indecomposable objects
$\ti{M}^k_n,\ti{M}^{k'}_{n'}\in\tHMF{R}$ and 
let us show that there exists a nonzero morphism 
$\Phi^\S\in\HomDz{R}(\ti{M}^{k}_{n},\S(\ti{M}^{k'}_{n'}))$ 
such that $\Phi^\S\Phi$ is nonzero in 
$\HomDz{R}(\ti{M}^{k'}_{n'},\S(\ti{M}^{k'}_{n'}))$. 
If $\ti{M}^{k}_{n}=\S(\ti{M}^{k'}_{n'})$, then 
one may take $\Phi^\S=\Id_{\S(\ti{M}^{k'}_{n'})}$, 
so assume that $\ti{M}^{k}_{n}\ne\S(\ti{M}^{k'}_{n'})$, 
\ie, $0\le \phi^{k}_{n}-\phi^{k'}_{n'}< 1-\frac{2}{h}$. 
As in the proof of Lemma \ref{lem:B}, 
from the long exact sequence (\ref{coh-functor}) of the 
AR-triangle (\ref{ARtri}) one obtains the exact sequence 
\begin{equation*}
0\to\HomDz{R}(\ti{M}^{k'}_{n'},\ti{M}^{k}_{n})
\to\HomDz{R}(\ti{M}^{k'}_{n'},C(\Psi))
\end{equation*} 
for the cone $C(\Psi)$ 
of a nonzero morphism $\Psi:\S^{-1}(\ti{M}^{k}_{n})\to\ti{M}^{k}_{n}$. 
Namely, there exists an indecomposable object 
$\ti{M}^{k''}_{n''}\in\tHMF{R}$ 
such that $d(k,k'')=1$, 
$\phi^{k''}_{n''}-\phi^{k}_{n}=\ov{h}$ and 
$\Phi'\Phi$ is not zero in 
$\HomDz{R}(\ti{M}^{k'}_{n'},\ti{M}^{k''}_{n''})$ 
with an irreducible morphism 
$\Phi'\in\HomDz{R}(\ti{M}^{k}_{n},\ti{M}^{k''}_{n''})$. 
Repeating this procedure 
together with Corollary \ref{cor:serre1} (a) leads that 
the path corresponding to the morphism arrives at 
$\S(\ti{M}^{k'}_{n'})$. 
Namely, for a given morphism $\Phi$ in 
$\HomDz{R}(\ti{M}^{k'}_{n'},\ti{M}^{k}_{n})$, 
there exists a nonzero morphism 
$\Phi^\S\in\HomDz{R}(\ti{M}^{k}_{n},\S(\ti{M}^{k'}_{n'}))$ 
such that $\Phi^\S\Phi$ is nonzero in 
$\HomDz{R}(\ti{M}^{k'}_{n'},\S(\ti{M}^{k'}_{n'}))$. 

Conversely, suppose that there exists a morphism 
$\Phi^\S\in\HomDz{R}(\ti{M}^{k}_{n},\S(\ti{M}^{k'}_{n'}))$ 
for two indecomposable objects 
$\ti{M}^{k}_{n},\ti{M}^{k'}_{n'}\in\tHMF{R}$. 
One can show that there exists a nonzero morphism 
$\Phi\in\HomDz{R}(\ti{M}^{k'}_{n'},\ti{M}^{k}_{n})$ 
such that $\Phi^\S\Phi$ is nonzero in 
$\HomDz{R}(\ti{M}^{k'}_{n'},\S(\ti{M}^{k'}_{n'}))$ by 
considering the functor $\HomDz{R}(-,\S(\ti{M}^{k'}_{n'}))$ 
on the AR-triangle (\ref{ARtri}) 
together with Corollary \ref{cor:serre1} (b), 
instead of $\HomDz{R}(\ti{M}^{k'}_{n'},-)$ 
with Corollary \ref{cor:serre1} (a). 
Thus, 
one gets Statement (iii-b) of Theorem \ref{thm:grading}. 
\qed\end{pf}

\vspace*{0.3cm}

\noindent
{\bf Step 3.\ \ Exceptional collections in $\tHMF{R}$}

In this step, for a fixed polynomial $f$ of type ADE 
and a Dynkin quiver $\vec{\Delta}$ 
of the corresponding type, 
we find an exceptional collection in $\tHMF{R}$ 
and then show the equivalence of the category $\tHMF{R}$ with 
the derived category of modules over 
the path algebra of $\vec{\Delta}$. 
\begin{lem}
For any Dynkin quiver $\vec{\Delta}$, 
there exists $\bn:=(n_1,\cdots,n_l)\in\Z^l$ such that 
one has an isomorphism of $\C$-algebras$:$ 
\begin{equation*}
\C\vec{\Delta}\simeq
\Abn:=\bigoplus_{k,k'\in\Pi} 
\HomDz{R}(\ti{M}^k_{n_k},\ti{M}^{k'}_{n_{k'}})\ ,
\end{equation*}
with the natural correspondence between 
the path from $k\in\Pi$ to $k'\in\Pi$ and 
the space $\Hom(\ti{M}^k_{n_k},\ti{M}^{k'}_{n_{k'}})$ 
for any $k,k'\in\Pi$. 
Such a tuple of integers is unique up to the $\Z$ shift 
$(n_1+n,\cdots,n_l+n)$ for some $n\in\Z$. 
 \label{lem:Calg}
\end{lem}
\begin{pf}
Let us fix $n_1\in\Z$ of $\ti{M}^1_{n_1}$. 
For each $k\in\Pi$, $d(1,k)=1$, 
take $\ti{M}^k_{n_k}\in\tHMF{R}$ 
such that $\phi^k_{n_k}-\phi^1_{n_1}=\ov{h}$ (resp. $-\ov{h}$) 
if there exists an arrow in $\vec{\Delta}$ from $1\in\Pi$ to $k\in\Pi$ 
(resp. from $k\in\Pi$ to $1\in\Pi$). 
Repeating this process leads to the tuple 
$\bn=(n_1,\cdots,n_l)$ 
such that 
$\Abn\simeq\C\vec{\Delta}$. 
On the other hand, 
Theorem \ref{thm:grading} (ii-b) implies that 
there exists 
a morphism for any given path and there are no relations 
among them. 
Thus, one gets this lemma. 
\qed\end{pf}
\begin{cor}
Let $\fS_l$ be the permutation group of\, $\Pi=\{1,\cdots,l\}$. 
Given a Dynkin quiver $\vec{\Delta}$ and 
$\bn\in\Z^l$ such that 
$\Abn\simeq\C\vec{\Delta}$, 
one can take an element $\sigma\in\fS_l$ such that 
$\phi^{\sigma(k')}_{n_{\sigma(k')}}\ge
\phi^{\sigma(k)}_{n_{\sigma(k)}}$ if 
$k'> k$. Thus, 
for
$\{E^1,\cdots,E^l\}$, $E^k:=\ti{M}^{\sigma(k)}_{n_{\sigma(k)}}$, 
$\HomDz{R}(E^k, E^{k'})\ne 0$ only if $k'>k$. 
 \label{cor:E}
\qed\end{cor}
\begin{lem}
Given a Dynkin quiver $\vec{\Delta}$ and 
$\bn\in\Z^l$ such that $\Abn\simeq\C\vec{\Delta}$, 
one has $\HomDz{R}(\tau^n(\ti{M}^k_{n_k}),\ti{M}^{k'}_{n_{k'}})=0$ 
for any $k,k'\in\Pi$ if $n\ge 1$. 
 \label{lem:hom-bound}
\qed\end{lem}
\begin{pf}
Let us first concentrate on irreducible morphisms. 
For any $k,k'\in\Pi$ such that $d(k,k')=1$, 
one obtains 
$\phi^{k'}_{n_{k'}}-\phi^k_{n_k+n}=\frac{\pm 1-2n}{h}$ 
since $\phi^{k'}_{n_{k'}}-\phi^k_{n_k}=\pm\ov{h}$. 
On the other hand, 
Theorem \ref{thm:grading} (ii-a) implies that there exists an irreducible 
morphism in $\HomDz{R}(\ti{M}^k_{n_k+n},\ti{M}^{k'}_{n_{k'}})$ 
only if $d(k,k')=1$ and $n=0$ or $n=-1$, 
since $\frac{\pm 1-2n}{h}$ can be $\ov{h}$ only if $n=0$ or $n=-1$. 

Since, by definition, all morphisms except the identities 
can be described by compositions of irreducible morphisms, 
one gets this lemma. 
\qed\end{pf}
By Lemma \ref{lem:hom-bound} and 
the Serre duality (Theorem \ref{thm:grading} (iii)), 
we obtain the following lemma.  
\begin{lem}
Given a Dynkin quiver $\vec{\Delta}$ and $\bn\in\Z^l$ such that 
$\Abn\simeq\C\vec{\Delta}$, 
$\{E^1,\dots, E^l\}$ given in Corollary \ref{cor:E} 
is a strongly exceptional collection, that is, 
\begin{equation*}
 \begin{cases}
 \HomDz{R}(E^i,E^j)=0 & \text{for}\ \, i>j \ ,\\ 
 \HomDz{R}(E^i,T^k(E^j))=0  &\text{for}\ \, k\ne 0\ \text{and any}\ \,i,j\ . 
 \end{cases}
\end{equation*}
 \label{lem:SEC}
\end{lem}
\begin{pf}
Due to Corollary \ref{cor:E}, 
it is sufficient to show that 
$\HomDz{R}(E^i,T^n(E^j))=0$ for $n\ge 1$. 
On the other hand, by the Serre duality, 
\begin{equation*}
\dim_\C(\HomDz{R}(E^i,T^n(E^j)))
 =\dim_\C(\HomDz{R}(T^{n-1}\tau(E^j),E^i))
\end{equation*}
holds. If $n=1$, 
$\HomDz{R}(T^{n-1}\tau(E^j),E^i)=\HomDz{R}(\tau(E^j),E^i)=0$ 
follows from Lemma \ref{lem:hom-bound}. 
If $n\ge 2$, 
first we have 
$\phi(T^{n-1}\tau(E^j))=\phi(E^j)+(n-1+\frac{2}{h})$ and then 
\begin{equation*}
 \phi(E^i)-\phi(T^{n-1}\tau(E^j))=(\phi(E^i)-\phi(E^j))
 -\left(n-1+\frac{2}{h}\right)\ .
\end{equation*}
Here, it is clear that 
\begin{equation*}
 |\phi(E^i)-\phi(E^j)|\le \frac{l-1}{h}
\end{equation*} 
holds for $\{E^1,\cdots,E^l\}$. 
Thus, we have the following inequality: 
\begin{equation*}
 \phi(E^i)-\phi(T^{n-1}\tau(E^j)) 
\le \frac{l-1}{h}-\left(n-1+\frac{2}{h}\right)< 0\quad \text{if}\ n\ge 2\ .
\end{equation*}
On the other hand, 
since all morphisms except the identities 
can be obtained by the composition of
the irreducible morphisms, Theorem \ref{thm:grading} (ii-a)
implies that 
$\HomDz{R}(E^j,E^i)\ne 0$ only if 
$\phi(E^i)-\phi(E^j)\ge\ov{h}d(E^j,E^i)\ge 0$. 
This implies \\ $\HomDz{R}(T^{n-1}\tau(E^j),E^i)=0$. 
\qed\end{pf}
\begin{cor}
$D^b(\mod\C\vec{\Delta})$ is a full triangulated 
subcategory of $\tHMF{R}$.
\end{cor}
\begin{pf}
Use the fact that $\{E^1,\dots,E^l\}$ is a strongly exceptional 
collection and 
$$
\C\vec{\Delta}\simeq \bigoplus_{i,j=1}^{l}\HomDz{R}(E^i,E^j)\ .
$$
Since $\tHMF{R}$ is an enhanced triangulated category, 
we can apply the theorem by Bondal-Kapranov $($\cite{bk:1} Theorem
1$)$ which implies that $D^b(\mod\C\vec{\Delta})$ is 
full as a triangulated subcategory. 
\qed\end{pf}
Recall that the number of indecomposable objects of $\tHMF{R}$ 
up to the shift functor $T$ 
is $\frac{l\cdot h}{2}$, 
which coincides with the number of the positive roots for 
the root system of type ADE.
This number is the number of indecomposable objects 
of $D^b(\mod\C\vec{\Delta})$ 
up to the shift functor $T$ 
by Gabriel's theorem \cite{g:1}. 
Therefore, one obtains $D^b(\mod\C\vec{\Delta})\simeq \tHMF{R}$. 

%%%%%%%%%%%%%%%%%%%%%%%%%%%%%%%%%%%%%%%%%%%%%%%%%%%%%%%%%%%%%%%%%%%%%%%%%%%%%%

\section{A stability condition on $\tHMF{R}$}
\label{sec:stab}
Let $K_0(\tHMF{R})$ be the Grothendieck group 
of the category $\tHMF{R}$ \cite{gr} (see also \cite{h}). 
For a triangulated category $\CC$, let $F$ be 
a free abelian group generated by 
the isomorphism classes of objects $[X]$ for $X\in\Ob(\CC)$. 
Let $F_0$ be the subgroup generated by 
$[X]-[Y]+[Z]$ for all  distinguished triangles 
$X\to Y\to Z\to T(X)$ in $\CC$. 
The Grothendieck group $K_0(\CC)$ 
of the triangulated category $\CC$ is, by definition, 
the factor group $F/F_0$. 
\begin{defn}\label{defn:centralcharge}
For a graded matrix factorization $\ti{M}=(Q,S)\in\tHMF{R}$, 
we associate a complex number as follows: 
\begin{equation*}
\ZZ(\ti{M}):= \Tr(e^{\pi\sqrt{-1}S})\ .
\end{equation*}
\end{defn}
The map $\ZZ$ extends to a group homomorphism 
$\ZZ: K_0(\tHMF{R})\to\C$. 
\begin{thm}\label{thm:stab}
Let $f\in R$ be a polynomial of type ADE. 
Let $\PP(\phi)$ be the full additive subcategory of 
$\tHMF{R}$ 
additively generated by the indecomposable objects of phase $\phi$. 
Then $\PP(\phi)$ and $\ZZ$ define a stability condition on 
$\tHMF{R}$ in the sense of Bridgeland \cite{bd:1}.
\end{thm}
\begin{pf}
By definition, what we should check is that 
$\PP(\phi)$ and $\ZZ$ satisfy the following properties$:$
\begin{enumerate}
\item if $\ti{M}\in \PP(\phi)$, then 
$\ZZ(\ti{M})=m(\ti{M})\exp(\pi\sqrt{-1}\phi)$ 
for some $m(\ti{M})\in\R_{> 0}$,
\item for all $\phi\in\R$, $\PP(\phi+1)=T(\PP(\phi))$,
\item if $\phi_1>\phi_2$ and $\ti{M}_i\in
\PP(\phi_i)$, $i=1,2$, then 
$\HomDz{R}(\ti{M}_1,\ti{M}_2)=0$,
\item for each nonzero object $\ti{M}\in \tHMF{R}$, 
there is a finite sequence of real numbers 
$$
\phi_1 >\phi_2>\dots >\phi_n
$$
and a collection of distinguished triangles
\[
\xymatrix@C=.5em{ \underset{}{\overset{}{0}} 
 \ar@{=}[r] & \underset{}{\overset{}{\ti{M}_0}} 
 \ar[rr] && \underset{}{\overset{}{\ti{M}_1}} \ar[rr]
\ar[dl] && \underset{}{\overset{}{\ti{M}_2}} 
 \ar[rr] \ar[dl] &&\ \ \ldots\ldots\ \ 
 \ar[rr] && \underset{}{\overset{}{\ti{M}_{n-1}}}
 \ar[rr] && \underset{}{\overset{}{\ti{M}_n}} 
 \ar[dl] \ar@{=}[r] &  \underset{}{\overset{}{\ti{M}}} \\
&& \ti{N}_1 \ar@{-->}[ul] && \ti{N}_2 \ar@{-->}[ul] &&&&&& \ti{N}_n
\ar@{-->}[ul] }
\]
with $\ti{N}_j\in \PP(\phi_j)$ for all $j=1,\cdots, n$. 
\end{enumerate}
Here, for any indecomposable object $\ti{M}^k_n\in\tHMF{R}$, 
$\ZZ(\ti{M}^k_n)$ is of the form: 
\begin{equation*}
 \ZZ(\ti{M}^k_n) =m(\ti{M}^k_n)e^{\pi\i\phi^k_n}\ ,\qquad 
 m(\ti{M}^k_n)\in\R_{> 0}\ .
\end{equation*}
In fact, by a straightforward calculation, one obtains 
\begin{equation*}
 \begin{split}
 \ZZ(\ti{M}^k_n)
 & =\sum_{j=1}^{\nu_k} 
 \left(e^{\pi\i(q^k_j+\phi^k_n)}+
 e^{\pi\i(-q^k_j+\phi^k_n)}+
 e^{\pi\i(\qb^k_j+\phi^k_n)}+
 e^{\pi\i(-\qb^k_j+\phi^k_n)} 
 \right) \\
 &=2 e^{\pi\i\phi^k_n}\sum_{j=1}^{\nu_k} 
 \left(\cos({\pi q^k_j})+  \cos({\pi\qb^k_j})\right)\\
 &=m(\ti{M}^k_n)e^{\pi\i\phi^k_n}\ ,\qquad 
 m(\ti{M}^k_n):=
 2\sum_{j=1}^{\nu_k}\left(\cos({\pi q^k_j})+  \cos({\pi\qb^k_j})\right)\ .
 \end{split}
\end{equation*}
This shows Statement (i) in Theorem \ref{thm:stab}. 
It is clear that 
Statement (ii),(iii) and (iv) follow from Theorem \ref{thm:grading}. 
\qed\end{pf}

Proposition 5.3 in \cite{bd:1} states that 
a stability condition $(\ZZ,\PP)$ 
on a triangulated category defines an abelian category. 
In our case, for the triangulate category $\tHMF{R}$, 
we obtain an abelian category $\PP((0,1])$ 
which consists of objects $\ti{M}\in\tHMF{R}$ 
having the decomposition $\ti{N}_1,\cdots,\ti{N}_n$ given by 
Theorem \ref{thm:stab} (iv) such that 
$1\ge\phi_1>\cdots>\phi_n>0$. 
\begin{prop}\label{prop:abelian}
Given a polynomial $f$ of type ADE, 
the following equivalence of abelian categories holds: 
\begin{equation*}
  \PP((0,1])\simeq \mod\C\vec{\Delta}_{principal}\ ,
\end{equation*}
where $\vec{\Delta}_{principal}$ is the Dynkin quiver of the
corresponding ADE Dynkin diagram $\Delta$ with the 
principal orientation \cite{sa:principal}, \ie, 
all arrows start from $\Pi_1$ and end on $\Pi_2$. 
\qed\end{prop}
\begin{pf}
Let us take the following collection 
$\{\ti{M}^1_0,\cdots,\ti{M}^l_0\}$. The corresponding $\C$-algebra 
is denoted by $\C\vec{\Gamma}(\bn=0)$ (see Lemma \ref{lem:Calg}). 
First, we show that
\begin{equation}\label{Calg-principal}
 \C\vec{\Gamma}(\bn=0)\simeq\C\vec{\Delta}_{principal}\ ,
\end{equation}
for which 
it is enough to say that $\ti{M}^k_0\in\PP((0,1])$ is projective 
for each $k\in\Pi$. 
For any $\ti{N}\in\PP(\phi(\ti{N}))$ such that $0<\phi(\ti{N})\le 1$, 
the Serre duality implies 
\begin{equation*}
 \HomDz{R}(\ti{M}^k_0,T(\ti{N}))
 \simeq \HomDz{R}(\tau(\ti{N}),\ti{M}^k_0)\ .
\end{equation*}
Here, one has $\HomDz{R}(\tau(\ti{N}),\ti{M}^k_0)\ne 0$ only if 
$\phi(\ti{N})+\frac{2}{h}\le\phi(\ti{M}^k_0)$. 
However, such $\ti{N}$ does not exist since 
$\phi(\ti{M}^k_0)=\frac{\sigma}{h}$ for $k\in\Pi_\sigma$, $\sigma=1,2$. 
Thus, one has $\HomDz{R}(\ti{M}^k_0,T(\ti{N}))=0$ and hence 
eq.(\ref{Calg-principal}). 

On the other hand, the full abelian subcategory 
$\langle\ti{M}^1_0,\cdots,\ti{M}^l_0\rangle$ of $\PP((0,1])$ 
generated by $\ti{M}^1_0,\cdots,\ti{M}^l_0$ 
is equivalent to the abelian category 
$\mod\C\vec{\Delta}_{principal}$. 
Also, the Gabriel's theorem \cite{g:1} asserts that 
the number of the indecomposable objects 
in $\mod\C\vec{\Delta}_{principal}$ 
is equal to the number of the positive roots of $\Delta$, 
which coincides with the number of the indecomposable 
objects in $\PP((0,1])$. 
Thus, one obtains the equivalence 
$\mod\C\vec{\Delta}_{principal}\simeq\PP((0,1])$. 
\qed\end{pf}
\begin{rem}[Principal orientation]
In the triangulated category $\tHMF{R}$, 
the principal oriented quiver $\vec{\Delta}_{principal}$ 
is realized by a strongly exceptional collection 
$\{\ti{M}^1_{n_1},\cdots,\ti{M}^l_{n_l}\}$ with 
$n_1=\cdots =n_l=n\in\Z$ for some $n\in\Z$ as above. 
It is interesting that 
a strongly exceptional collection of this type 
has minimum range of the phase: 
$\frac{2n+1}{h}\le\phi(\ti{M}^k_{n_k})\le\frac{2n+2}{h}$ 
for any $k\in\Pi$. 
 \label{rem:principal}
\end{rem}

 \section{Tables of data for matrix factorizations of type ADE}
\label{sec:tbl}

\begin{tbl}
We list up the Auslander-Reiten (AR)-quiver 
of the category $\HMF{\Oh}$ for a polynomial 
$f$ of type ADE. 
We label each isomorphism class of the indecomposable objects 
(matrix factorizations) in $\MF{\Oh}$ 
by upperscript $\{1,\cdot\cdot,k,\cdot\cdot,l\}$ such as $[M^k]$. 
We also list up a representative $M^k$ of $[M^k]$ 
by giving the pair $(\varphi^k,\psi^k)$ for each 
matrix factorization $Q^k:=\left(\bps  & \psi^k\\ \varphi^k & \eps\right)$ 
(see eq.(\ref{Q})). 
In addition, for type $D_l$ and $E_l\ (l=6,7,8)$ cases, 
we attach the indices 
$\left(\bps q^k_1\\ \vdots \\ q^k_{\nu_k}\eps\right)$ 
defined in Theorem \ref{thm:grading} (eq.(\ref{gradADE})) 
to each of those pairs as 
\begin{equation*}
 \varphi^k_{\left(\bps q^k_1\\ \vdots \\ q^k_{\nu_k}\eps\right)},
\psi^k_{\left(\bps q^k_1\\ \vdots \\ q^k_{\nu_k}\eps\right)} \ .
\end{equation*}
Those indices are for later convenience and 
irrelevant for the statement of Theorem \ref{thm:AR}. 
They will be listed up again in the next table 
(Table \ref{tbl:2}).

\noindent
{\bf Type $A_l$.}\ \, The AR-quiver for the polynomial 
$f=x^{l+1}+yz$ is of the form: 
\begin{equation*}
\xymatrix{
 \ \ [M^1]\ \ar@<0.5ex>[r]^{}
  & \ [M^2]\ \ar@<0.5ex>[l]^{}\ar@<0.5ex>[r]^{}
  & \ \cdot\cdots\cdot\ \ar@<0.5ex>[l]^{}\ar@<0.5ex>[r]^{}
  & \ [M^{l-1}]\ \ar@<0.5ex>[l]^{}\ar@<0.5ex>[r]^{}
  & \ [M^l]\ \ar@<0.5ex>[l]^{}\ \ ,
}
\end{equation*}
where, for each $k\in\Pi$, $(\varphi^k,\psi^k)$ is given by 
{\small 
\begin{equation*}
\varphi^k:=
\begin{pmatrix}
y & x^{l+1-k}\\
x^k & -z
\end{pmatrix},\quad \psi^k:=
\begin{pmatrix}
z & x^{l+1-k}\\
x^k & -y
\end{pmatrix}\ .
\end{equation*}
}
\vspace*{0.3cm}

\noindent
{\bf Type $D_l$.}\ \, The AR-quiver for the polynomial
$f=x^2y+y^{l-1}+z^2$ is given by: 
\begin{equation*}
\xymatrix@C=.6em{
 && && && && &\ [M^{l-1}]\  \ar@<0.5ex>[ld]^{}\\
\ [M^1]\ \ar@<0.5ex>[rr]^{} 
  && \ [M^2]\ \ar@<0.5ex>[ll]^{}
 \ar@<0.5ex>[rr]^{} 
  && \ \  {\cdots} \ \  \ar@<0.5ex>[ll]^{}  
  \cdots \ar@<0.5ex>[rr]^{}  
  &&\ [M^{l-3}]\ 
 \ar@<0.5ex>[rr]^{} \ar@<0.5ex>[ll]^{} 
  &&\ [M^{l-2}]\ 
 \ar@<0.5ex>[ru]^{} \ar@<0.5ex>[rd]^{}
 \ar@<0.5ex>[ll]^{} \\
 && && && && & 
 \ \ \ [M^l]\ \ar@<0.5ex>[lu]^{}  \ \ ,
}
\end{equation*}
where, for each 
$k\in\Pi$, $(\varphi^k,\psi^k)$ is given by 
{\small 
\begin{equation*}
  \varphi^1_{\left(\bps l-3\eps\right)}=
  \psi^1_{\left(\bps l-3\eps\right)}=
\bp z & x^2+y^{l-2}\\
    y & -z \ep\ ,
\end{equation*}
\begin{equation*}
  \varphi^k_{\left(\bps l-2-k\\l-k \eps\right)}=
  \psi^k_{\left(\bps l-2-k\\l-k \eps\right)}=
 \bp    -z &   0 &        xy & y^{\frac{k}{2}} \\
         0 &  -z & y^{l-1-\frac{k}{2}} & -x \\
         x & y^\frac{k}{2} &         z & 0 \\ 
 y^{l-1-\frac{k}{2}} & -xy &         0 & z \ep\ ,\quad 
 k \text{:\ even}\ \ (2\le k\le l-2)
\end{equation*}
\begin{equation*}
  \varphi^k_{\left(\bps l-2-k\\l-k \eps\right)}=
  \psi^k_{\left(\bps l-2-k\\l-k \eps\right)}=
 \bp -z      & y^{\frac{k-1}{2}} &        xy & 0 \\ 
   y^{l-\frac{k+1}{2}} &       z &         0 & -x \\
           x &       0 &         z & y^\frac{k-3}{2} \\ 
           0 &     -xy & y^{l-\frac{k-1}{2}} & -z \ep\ ,\quad 
 k \text{:\ odd}\ \ (3\le k\le l-2)
\end{equation*}
}The forms of $(\varphi^{l-1},\psi^{l-1})$ and $(\varphi^l,\psi^l)$ 
depend on whether $l$ is even or odd. 

\noindent
$\bullet$\ If $l$ is even, one obtains 
{\small 
\begin{equation*}
 \begin{split}
 & \varphi^{l-1}_{\left(\bps 1 \eps\right)}=
  \psi^{l-1}_{\left(\bps 1 \eps\right)}=
 \bp z                        & y(x+\i y^{\frac{l-2}{2}}) \\ 
     x-\i y^{\frac{l-2}{2}} & -z 
 \ep\ ,\\
 & \varphi^l_{\left(\bps 1 \eps\right)}=
  \psi^l_{\left(\bps 1 \eps\right)}=
 \bp z                        & y(x-\i y^{\frac{l-2}{2}}) \\ 
     x+\i y^{\frac{l-2}{2}} & -z 
 \ep\ .
 \end{split}
\end{equation*}
}
\noindent
$\bullet$\ If $l$ is odd, one obtains 
{\small 
\begin{equation*}
 \begin{split}
  &\varphi^{l-1}_{\left(\bps 1 \eps\right)}
  =\psi^{l}_{\left(\bps 1 \eps\right)}=
\bp z+\i y^{\frac{l-1}{2}} & xy \\
                           x & -(z-\i y^{\frac{l-1}{2}})
 \ep\ ,\\
  &\varphi^l_{\left(\bps 1 \eps\right)}
  =\psi^{l-1}_{\left(\bps 1 \eps\right)}=
\bp z-\i y^{\frac{l-1}{2}} & xy \\
                           x & -(z+\i y^{\frac{l-1}{2}})
 \ep\ . 
 \end{split}
\end{equation*}
}

\vspace*{0.3cm}

\noindent
{\bf Type $E_6$.}\ \, The AR-quiver for the polynomial 
$f=x^3+y^4+z^2$ is given by: 
\begin{equation*}
\xymatrix{
 & & \ [M^1] 
 \ar@<0.7ex>[d]^{} & \\
 \ [M^5]\ \ar@<0.5ex>[r]^{}
 & \ [M^3] \ \ar@<0.5ex>[r]^{}
 \ar@<0.5ex>[l]^{}
 & \ [M^2] \ \ar@<0.5ex>[r]^{}
 \ar@<0.5ex>[l]^{} \ar@<0.5ex>[u]^{}
 & \ [M^4] \ \ar@<0.5ex>[r]^{}
 \ar@<0.5ex>[l]^{}
 & \ [M^6] \ \ . \ar@<0.5ex>[l]^{}
}
\end{equation*}
For $Y^\pm:=y^2\pm\i z$, each pair $(\varphi^k,\psi^k)$ is given by 
{\footnotesize 
\begin{equation*}
 \varphi^1_{\left(\bps 1 \\ 5\eps\right)}
 =\varphi^1_{\left(\bps 1 \\ 5\eps\right)}
 =\bp 
 -z & 0 & x^2 & y^3 \\ 
 0 & -z & y & -x \\
 x & y^3 & z & 0  \\
 y & -x^2 & 0 & z 
 \ep\ , 
\end{equation*}
\begin{equation*}
 \varphi^2_{\left(\bps 0\\ 2\\ 4\eps\right)}=
 \bp -\i z &  -y^2 &    xy &     0 &   x^2 & 0 \\ 
      -y^2 & -\i z &     0 &     0 &     0 & x \\
         0 &     0 & -\i z &    -x &     0 & y \\
         0 &    xy &  -x^2 & -\i z &   y^3 & 0 \\
         x &     0 &     0 &     y & -\i z & 0 \\
         0 &   x^2 &   y^3 &     0 &  xy^2 & -\i z
\ep\ ,
\end{equation*}
\begin{equation*}
\psi^2_{\left(\bps 0\\ 2\\ 4\eps\right)}=
 \bp  \i z &  -y^2 &    xy &     0 &   x^2 & 0 \\ 
      -y^2 &  \i z &     0 &     0 &     0 & x \\
         0 &     0 &  \i z &    -x &     0 & y \\
         0 &    xy &  -x^2 &  \i z &   y^3 & 0 \\
         x &     0 &     0 &     y &  \i z & 0 \\
         0 &   x^2 &   y^3 &     0 &  xy^2 & \i z
 \ep\ ,
\end{equation*}
\begin{equation*}
\varphi^3_{\left(\bps 1\\ 3\eps\right)}=
 \psi^4_{\left(\bps 1\\ 3\eps\right)}=
 \bp -Y^- &     0 &      xy & x\\ 
        -xy & Y^+ &     x^2 & 0 \\
          0 &     x & \i z & y \\
        x^2 &   -xy &     y^3 & \i z \ep\ ,
\end{equation*}
\begin{equation*}
\varphi^4_{\left(\bps 1\\ 3\eps\right)}=
 \psi^3_{\left(\bps 1\\ 3\eps\right)}=
 \bp -Y^+ &     0 &      xy & x\\ 
        -xy & Y^- &     x^2 & 0 \\
          0 &     x & -\i z & y \\
        x^2 &   -xy &     y^3 & -\i z \ep\ ,
\end{equation*}
\begin{equation*}
\varphi^5_{\left(\bps 2 \eps\right)}= 
  \psi^6_{\left(\bps 2 \eps\right)}= 
  \bp  -Y^- & x      \\
          x^2 & Y^+ \ep\ ,\qquad 
\varphi^6_{\left(\bps 2 \eps\right)}= 
  \psi^5_{\left(\bps 2 \eps\right)}= 
  \bp  -Y^+ & x      \\
          x^2 & Y^- \ep\ . 
\end{equation*}
}

\vspace*{0.3cm}

\noindent
{\bf Type $E_7$.}\ \, The AR-quiver for the polynomial 
$f=x^3+xy^3+z^2$ is given by: 
\begin{equation*}
\xymatrix{
 & & & [M^4] \ar@<0.5ex>[d]^{} & & & \\
 \ [M^7] \  \ar@<0.5ex>[r]^{} 
 & \ [M^6] \  \ar@<0.5ex>[r]^{} 
 \ar@<0.5ex>[l]^{} 
 & \ [M^5] \  
 \ar@<0.5ex>[r]^{} \ar@<0.5ex>[l]^{} 
 & \ [M^3] \  \ar@<0.5ex>[r]^{} 
 \ar@<0.5ex>[l]^{} \ar@<0.5ex>[u]^{}
 & \ [M^2] \  
\ar@<0.5ex>[r]^{} \ar@<0.5ex>[l]^{} 
 & \ [M^1] \ . 
\ar@<0.5ex>[l]^{} 
 & 
}
\end{equation*}
The corresponding matrix factorizations are: 
{\footnotesize 
\begin{equation*}
 \varphi^1_{\left(\bps 2\\8\eps \right)}
 = \psi^1_{\left(\bps 2\\8\eps \right)}
 =\bp z & 0 & -x^2 & y \\
      0 & z & xy^2 & x \\ 
      -x & y & -z & 0 \\
      xy^2 & x^2 & 0 & -z \ep\ ,\quad
 \varphi^4_{\left(\bps 1\\5\eps \right)}
 = \psi^4_{\left(\bps 1\\5\eps \right)}
 =\bp -z & y^2 & 0 & x \\
      xy & z & -x^2 & 0\\ 
      0 & -x & -z & y \\
      x^2 & 0 & xy^2 & z \ep\ 
\end{equation*}
\begin{equation*}
 \varphi^2_{\left(\bps 1\\3\\7\eps \right)}
 = \psi^2_{\left(\bps 1\\3\\7\eps \right)}
 =\bp -z & y^2 & xy & 0 & x^2 & 0 \\
      xy & z & 0 & 0 & 0 & -x\\ 
      0 & 0 & z & -x & 0 & y \\
      0 & -xy & -x^2 & -z & xy^2 & 0 \\
      x & 0 & 0 & y & z & 0 \\
      0 & -x^2 & xy^2 & 0 & x^2y & -z \ep\ ,
\end{equation*}
\begin{equation*}
 \varphi^5_{\left(\bps 1\\3\\5\eps \right)}
 =\psi^5_{\left(\bps 1\\3\\5\eps \right)}
 =\bp -z & 0 & xy & 0 & 0 & x \\
      -xy & z & 0 & -y^2 & -x^2 & 0\\ 
      y^2 & 0 & z & -x & xy & 0 \\
      0 & -xy & -x^2 & -z & 0 & 0 \\
      0 & -x & 0 & 0 & -z & -y \\
      x^2 & 0 & 0 & xy & -xy^2 & z \ep\ ,
\end{equation*}
\begin{equation*}
 \varphi^3_{\left(\bps 0\\2\\4\\6\eps \right)}
 = \psi^3_{\left(\bps 0\\2\\4\\6\eps \right)}
 =\bp 
 -z & 0 & xy & -y^2 & 0 & 0 & x^2 & 0 \\
 0 & -z & 0 & y^2 & 0 & 0 & 0 & x \\
 y^2 & y^2 & z & 0 & 0 & -x & 0 & 0\\
 0 & xy & 0 & z & -x^2 & 0 & 0 & 0\\
 0 & 0 & 0 & -x & -z & 0 & 0 & y \\
 0 & 0 & -x^2 & 0 & 0 & -z & xy^2 & y^2 \\
 x & 0 & 0 & 0 & -y^2 & y^2 & z & 0\\
 0 & x^2 & 0 & 0 & xy^2 & 0 & 0 & z
\ep\ ,
\end{equation*}
\begin{equation*}
 \varphi^6_{\left(\bps 2\\4\eps \right)}
 = \psi^6_{\left(\bps 2\\4\eps \right)}
 =\bp z & 0 & -xy & x \\
      0 & z & x^2 & y^2\\ 
      -y^2 & x & -z & 0 \\
      x^2 & xy & 0 & -z \ep\ ,\quad 
 \varphi^7_{\left(\bps 3\eps \right)}
 = \psi^7_{\left(\bps 3\eps \right)}
 =\bp z & x \\
      x^2+y^3 & -z \ep\ .
\end{equation*}
}

\vspace*{0.3cm}

\noindent
{\bf Type $E_8$.}\ \, 
The AR-quiver for the polynomial $f=x^3+y^5+z^2$ is given by: 
\begin{equation*}
\hspace*{-0.7cm}
\xymatrix{
 & & & & [M^6] \ar@<0.5ex>[d]^{} & & \\
 \ [M^1] \  \ar@<0.5ex>[r]^{} 
 & \ [M^2] \  
 \ar@<0.5ex>[r]^{} \ar@<0.5ex>[l]^{} 
 & \ [M^3] \  
 \ar@<0.5ex>[r]^{} \ar@<0.5ex>[l]^{} 
 & \ [M^4] \  
\ar@<0.5ex>[r]^{} \ar@<0.5ex>[l]^{} 
 & \ [M^5] \  
\ar@<0.5ex>[r]^{} \ar@<0.5ex>[l]^{} \ar@<0.5ex>[u]^{}
 & \ [M^7] \  
\ar@<0.5ex>[r]^{} \ar@<0.5ex>[l]^{} 
 & \ [M^8] \ . 
\ar@<0.5ex>[l]^{} 
}
\end{equation*}
The corresponding matrix factorizations are : 
{\scriptsize 
\begin{equation*}
 \hspace*{-0.5cm}
 \varphi^1_{\left(\bps 4\\14\eps \right)}
 = \psi^1_{\left(\bps 4\\14\eps \right)}
 =\bp z & 0 & x & y \\
      0 & z & y^4 & -x^2\\ 
      x^2 & y & -z & 0 \\
      y^4 & -x & 0 & -z \ep\ ,\quad 
 \varphi^8_{\left(\bps 2\\8\eps \right)}
 = \psi^8_{\left(\bps 2\\8\eps \right)}
 =\bp z & 0 & x & y^2 \\
      0 & z & y^3 & -x^2\\ 
      x^2 & y^2 & -z & 0 \\
      y^3 & -x & 0 & -z \ep\ ,
\end{equation*}
\begin{equation*}
 \varphi^2_{\left(\bps 3\\5\\13\eps \right)}
 = \psi^2_{\left(\bps 3\\5\\13\eps \right)}
 =\bp 
 z & -y^2 & xy & 0 & -x^2 & 0 \\
 -y^3 & -z & 0 & 0 & 0 & x \\
 0 & 0 & -z & x & 0 & y \\
 0 & -xy & x^2 & z & y^4 & 0\\
 -x & 0 & 0 & y & -z & 0 \\
 0 & x^2 & y^4 & 0 & -xy^3 & z
 \ep\ ,
 \ \ 
 \varphi^6_{\left(\bps 1\\5\\9\eps \right)}
 =\psi^6_{\left(\bps 1\\5\\9\eps \right)}
 =\bp 
 -z & 0 & 0 & y^2 & 0 & x \\
 xy & z & -y^3 & 0 & -x^2 & 0 \\
 0 & -y^2 & -z & x & 0 & 0 \\
 y^3 & 0 & x^2 & z & -xy^2 & 0\\
 0 & -x & 0 & 0 & -z & y \\
 x^2 & 0 & -xy^2 & 0 & y^4 & z
 \ep\ ,
\end{equation*}
\begin{equation*}
 \varphi^3_{\left(\bps 2\\4\\6\\12\eps \right)}
 = \psi^3_{\left(\bps 2\\4\\6\\12\eps \right)}
 =\bp 
 -z & 0 & -xy & y^2 & 0 & 0 & x^2 & 0 \\
 0 & -z & y^3 & 0 & 0 & 0 & 0 & x \\
 0 & y^2 & z & 0 & 0 & -x & 0 & 0\\
 y^3 & xy & 0 & z & -x^2 & 0 & 0 & 0\\
 0 & 0 & 0 & -x & -z & 0 & y^3 & y \\
 0 & 0 & -x^2 & 0 & 0 & -z & 0 & y^2 \\
 x & 0 & 0 & 0 & y^2 & -y & z & 0\\
 0 & x^2 & 0 & 0 & 0 & y^3 & 0 & z
\ep\ ,
\end{equation*}
\begin{equation*}
 \varphi^4_{\left(\bps 1\\3\\5\\7\\11\eps \right)}
= \psi^4_{\left(\bps 1\\3\\5\\7\\11\eps \right)}
 =\bp 
 z & 0 & xy & 0 & 0 & -y^2 & y^3 & 0 & -x^2 & 0 \\
 0 & -z & 0 & 0 & 0 & 0 & 0 & -y^2 & 0 & x \\
 0 & 0 & -z & y^2 & 0 & 0 & 0 & x & 0 & 0\\
 0 & xy & y^3 & z & 0 & 0 & -x^2 & 0 & 0 & 0\\
 0 & y^2 & 0 & 0 & z & -x & 0 & 0 & y^3 & 0 \\
 -y^3 & 0 & 0 & 0 & -x^2 & -z & 0 & 0 & 0 & y^2 \\
 0 & 0 & 0 & -x & 0 & 0 & -z & 0 & 0 & y\\
 0 & -y^3 & x^2 & 0 & 0 & 0 & xy^2 & z & 0 & 0 \\
 -x & 0 & 0 & 0 & y^2 & 0 & 0 & y & -z & 0\\
 0 & x^2 & xy^2 & 0 & 0 & 0 & y^4 & 0 & 0 & z 
\ep\ ,
\end{equation*}
\begin{equation*}
\hspace*{-0.8cm}
 \varphi^5_{\left(\bps 0\\2\\4\\6\\8\\10\eps \right)}
 = \psi^5_{\left(\bps 0\\2\\4\\6\\8\\10\eps \right)}
 =
 \left(
 \begin{array}{cccccccccccc}
 -z & 0 & 0 & 0 & 0 & 0 & 0 & y^2 & 0 & 0 & 0 & x \\
 0 & -z & -xy & 0 & 0 & 0 & y^3 & -y^2 & 0 & 0 & x^2 & 0 \\
 0 & 0 & z & 0 & 0 & -y^2 & 0 & 0 & y^3 & -x & 0 & 0 \\
 xy & 0 & 0 & z & -y^3 & 0 & 0 & 0 & -x^2 & 0 & 0 & 0 \\
 0 & 0 & 0 & -y^2 & -z & 0 & 0 & x & 0 & 0 & 0 & 0 \\
 0 & 0 & -y^3 & 0 & 0 & -z & -x^2 & 0 & 0 & 0 & xy^2 & y^2 \\
%6
 y^2 & y^2 & 0 & 0 & 0 & -x & z & 0 & 0 & 0 & 0 & 0 \\
 y^3 & 0 & 0 & 0 & x^2 & 0 & 0 & z & -xy^2 & 0 & 0 & 0 \\
%8
 0 & 0 & 0 & -x & 0 & 0 & 0 & 0 & -z & 0 & 0 & y \\
 0 & 0 & -x^2 & -y^3 & 0 & 0 & xy^2 & 0 & 0 & -z & -y^4 & 0 \\
%10
 0 & x & 0 & 0 & y^2 & 0 & 0 & 0 & 0 & -y & z & 0 \\
 x^2 & 0 & 0 & 0 & -xy^2 & 0 & 0 & 0 & y^4 & 0 & 0 & z 
 \end{array}\right).
\end{equation*}
\begin{equation*}
 \varphi^7_{\left(\bps 1\\3\\7\\9\eps \right)}
 =\psi^7_{\left(\bps 1\\3\\7\\9\eps \right)}
 =\bp 
 z & 0 & 0 & 0 & -y^3 & 0 & 0 & -x \\
 xy & -z & 0 & 0 & 0 & y^2 & x^2 & 0 \\
 0 & 0 & -z & y^2 & 0 & x & -y^3 & 0\\
 0 & 0 & 0 & z & -x^2 & 0 & 0 & y^2\\
 -y^2 & 0 & 0 & -x & -z & 0 & 0 & 0 \\
 0 & y^3 & x^2 & 0 & xy^2 & z & 0 & 0 \\
 0 & x & -y^2 & 0 & 0 & 0 & z & y\\
 -x^2 & 0 & 0 & y^3 & 0 & 0 & 0 & -z
\ep\ .
\end{equation*}
}
The shift functor $T$ acts on these matrix factorizations as follows. 
\begin{itemize}
 \item 
For type $A_l$, $T(M^k)\simeq M^{l+2-k}$ for any $k\in\Pi$. 
 \item 
For type $D_l$, $T(M^k)\simeq M^k$ for all $k\in\Pi$, 
except that $T(M^{l-1})\simeq M^l$, $T(M^l)\simeq M^{l-1}$ 
if $l$ is odd. 
 \item
For type $E_6$, $T(M^k)\simeq M^k$ for $k=1,2$ but $T(M^3)\simeq M^4$, 
$T(M^4)\simeq M^3$ and 
$T(M^5)\simeq M^6$, $T(M^6)\simeq M^5$.
 \item
For type $E_7$ or $E_8$, $T(M^k)\simeq M^k$ for any $k\in\Pi$. 
\end{itemize}
Remember that the Serre functor 
(in Theorem \ref{thm:grading} (iii)) 
is defined by $\S:=T\tau^{-1}$ when we shall introduce the grading. 
 \label{tbl:1}
\end{tbl}

\begin{tbl}
The list of all the isomorphism classes of the indecomposable objects 
in $\tHMF{R}$ for a polynomial $f\in R$ of type ADE is given. 

The set of the isomorphism classes of the indecomposable objects 
is given by 
\begin{equation*}
\left\{[\ti{M}^k_n:=(Q^k,S^k_n)],
\quad k\in\Pi,\quad n\in\Z \right\}\ . 
\end{equation*}
Here, for each $k\in\Pi$, $Q^k$ is the matrix factorization 
of size $2\nu_k$ given in Table \ref{tbl:1}, 
and the grading matrix $S^k_n$ is a diagonal matrix as follows: 
\begin{equation*}
 S^k_n:
=\diag\left(q^k_1,-q^k_1,\cdots,q^k_{\nu_k},-q^k_{\nu_k};
 \qb^k_1,-\qb^k_1,\cdots,\qb^k_{\nu_k},-\qb^k_{\nu_k}\right)
+\phi^k_n\cdot\1_{4\nu_k}\ ,
\end{equation*}
where the phase is given by 
$\phi^k_n=\frac{2n+\sigma}{h}$ for $k\in\Pi_\sigma$, $\sigma=1,2$, 
and the data $q^k_j,\qb^k_j$ are given below. 

 \vspace*{0.3cm}

\noindent
{\bf Type $A_l$}\ $(h=l+1)$: 
\quad 
In this case, $\nu_k=1$ for all $k\in\Pi$ 
and the grading is given by 
\begin{equation*}
(q^k_1;\qb^k_1)
 =\ov{l+1}(b-k;(l+1-b)-k)\ . 
\end{equation*}

 \vspace*{0.3cm}

For type $D_l$ and $E_l$, for any $k\in\Pi$ 
there exists a representative of 
the isomorphism classes of the indecomposable objects such that 
$q^k_j=\qb^k_j$, $j=1,\cdots,\nu_k$. 
The matrix factorizations $Q^k$'s listed in Table \ref{tbl:1} 
are just such ones. 
Therefore, we present only $q^k_j$ and omit $\qb^k_j$. 

 \vspace*{0.3cm}

\noindent
{\bf Type $D_l$}\ $(h=2(l-1))$: 
\begin{minipage}[t]{80mm}{
\begin{equation*}
 \begin{array}{lll}
 k & \nu_k & (q^k_1,\cdots,q^k_{\nu_k}) \\
 1 &   1 & \ov{2(l-1)}(l-3) \\
 2,\cdots,(l-2)& 2 & \ov{2(l-1)}(l-k-2,l-k) \\
 (l-1),l & 1 & \ov{2(l-1)}(1) \ .\\
 \end{array}
\end{equation*}
}
\end{minipage}
 \vspace*{0.3cm}

\noindent

\begin{equation*}
 \begin{array}{lll}
 \text{{\bf Type $E_6$}\ $(h=12)$: }\qquad\quad  
 & \text{{\bf Type $E_7$}\ $(h=18)$: }\qquad\quad
 & \text{{\bf Type $E_8$}\ $(h=30)$: }\qquad \\
 \begin{array}{llc}
 k & \nu_k & (q^k_1,\cdots,q^k_{\nu_k}) \\
 1 & 2 & \ov{12}(1,5)\ \\
 2 & 3 & \ov{12}(0,2,4)\ \\
 3,4 & 2 & \ov{12}(1,3)\ \\
 5,6 & 1 & \ov{12}(2)\ .
 \end{array}
 &
 \begin{array}{llc}
 k & \nu_k & (q^k_1,\cdots,q^k_{\nu_k}) \\
 1 & 2   & \ov{18}(2,8)\ \\
 2 & 3   & \ov{18}(1,3,7)\ \\
 3 & 4   & \ov{18}(0,2,4,6)\ \\
 4 & 2   & \ov{18}(1,5)\ \\
 5 & 3   & \ov{18}(1,3,5)\ \\
 6 & 2   & \ov{18}(2,4)\ \\
 7 & 1   & \ov{18}(3)\ .
 \end{array}
 &
 \begin{array}{ccc}
 k & \nu_k & (q^k_1,\cdots,q^k_{\nu_k}) \\
 1 & 2 & \ov{30}(4,14)\ \\
 2 & 3 & \ov{30}(3,5,13)\ \\
 3 & 4 & \ov{30}(2,4,6,12)\ \\
 4 & 5 & \ov{30}(1,3,5,7,11)\ \\
 5 & 6 & \ov{30}(0,2,4,6,8,10)\ \\
 6 & 3 & \ov{30}(1,5,9)\ \\
 7 & 4 & \ov{30}(1,3,7,9)\ \\
 8 & 2 & \ov{30}(2,8)\ .
 \end{array}
 \end{array}
\end{equation*}
\begin{enumerate}
 \item If we set $n=0$ for all $k\in\Pi$, then 
$\{\ti{M}^1_0,\cdots,\ti{M}^l_0\}$ corresponds to the Dynkin quiver
with a principal orientation (see Proposition \ref{prop:abelian}).

 \item For the grading matrix of $\ti{M}^k_n$, 
 one has $q^k_j\ne 0$ for $j>1$, and $q^k_1=0$ if
 and only if $F(\ti{M}^k_n)=M_o$. 

 \item Lemma \ref{lem:B} can be checked 
at the level of the grading matrices $S$. Namely, 
for an indecomposable object $\ti{M}^k_n\in\tHMF{R}$, 
the cone $C(\Psi)$ of a nonzero morphism 
$\Psi\in\HomDz{R}(\S^{-1}(\ti{M}^k_n),\ti{M}^k_n)$ 
is isomorphic to 
the direct sum of indecomposable objects 
$\ti{M}^{k_1}_{n_1}\oplus\cdots\oplus\ti{M}^{k_m}_{n_m}$ 
for some $m\in\Z_{>0}$ 
such that $\{k_1,\cdots,k_m\}=\{k'\in\Pi\ |\ d(k,k')=1\}$ and 
$\phi(\ti{M}^{k_i}_{n_i})=\phi(\ti{M}^k_n)+\ov{h}$ 
for any $i=1,\cdots,m$. 
Correspondingly, one can check that 
\begin{align}
 & \left\{q^k_1-\ov{h},-q^k_1-\ov{h},\cdots,q^k_{\nu_k}-\ov{h},
 -q^k_{\nu_k}-\ov{h}\right\}
\coprod 
\left\{q^k_1+\ov{h},-q^k_1+\ov{h},\cdots,q^k_{\nu_k}+\ov{h},
 -q^k_{\nu_k}+\ov{h}\right\} \nonumber\\
 & 
 \qquad=\coprod_{i=1}^m
 \left\{q^{k_i}_1,-q^{k_i}_1,\cdots,
 q^{k_i}_{\nu_{k_i}},-q^{k_i}_{\nu_{k_i}}\right\}  \label{cone-grading}
\end{align}
\end{enumerate}
 \label{tbl:2}
\end{tbl}
and the same identities for $\qb^k_j$, $k\in\Pi$, $j=1,\cdots,\nu_k$.

\begin{tbl}\label{tbl:3} 
The following tables give all the morphisms between 
all the isomorphism classes of the indecomposable objects in $\tHMF{R}$ 
for a polynomial $f$ of type ADE. 

Recall that, 
for two indecomposable objects $\ti{M}^k_n,\ti{M}^{k'}_{n'}\in\tHMF{R}$, 
$\fC(k,k')$ is the multi-set of non-negative integers such that 
\begin{equation*}
\fC(k,k')
 :=\{c:=h(\phi^{k'}_{n'}-\phi^k_n)+2n''\ |\ 
 \HomDz{R}(\ti{M}^k_n,\tau^{n''}(\ti{M}^{k'}_{n'}))\ne 0\ ,\ \ n\in\Z\}\ ,
\end{equation*} 
where the integer $c$ appears with multiplicity 
$d:=\dim_\C(\HomDz{R}(\ti{M}^k_n,\tau^{n''}(\ti{M}^{k'}_{n'})))$. 
We sometimes write $c^d$ instead of $d$ copies of $c$. 

For each type of ADE, $\fC(k,k')$ is listed up for any 
$k,k'\in\Pi$.

\vspace*{0.3cm}

\noindent 
{\bf Type $A_l$}\ $(h-2=l-1)$\ : \  
One obtains (essentially the same as those given in 
\cite{hw:1}) 
\begin{equation*}
 \fC(k,k')=\left\{
 \begin{split}
 & |k'-k|,\ |k'-k|+2,\ |k'-k|+4,\ \cdots\cdots\ , \\
 & \hspace*{1.5cm}\cdots\cdots\ ,l-3-|(l-1)-(k+k'-2)|,\ \, 
 l-1-|(l-1)-(k+k'-2)| 
 \end{split}
 \right\}\ .
\end{equation*}
For any $k,k'\in\Pi$, 
$\fC(k,k')$ does not contain multiple copies of the same integer. 
Pictorially, the table of $\fC(k,k')$ is displayed as 
\begin{equation*}
 \renewcommand{\arraystretch}{1.5}
\begin{array}{@{\vrule width 1pt\ \ }c|c|c|c|c|c|c|c@{\ \vrule width 1pt}}
  \hline
 \overset{}{k}\backslash \underset{}{k'} & 
 1 & 2 & 3 & \cdots\cdots\cdots & l-2 & l-1 & l \\ 
  \hline 
 1 & 
 \bps 0 \eps & 
 \bps 1 \eps & 
 \bps 2 \eps & 
 \bps \cdots\cdots\cdots \eps & 
 \bps l-3 \eps & 
 \bps l-2 \eps & 
 \bps l-1 \eps \\
  \hline
 2 & 
 \bps 1 \eps & 
 \bps 0 & 2 \eps & 
 \bps 1 & 3 \eps & 
 \bps \cdots\cdots\cdots \eps & 
 \bps l-4 & l-2 \eps &
 \bps l-3 & l-1 \eps &
 \bps l-2 \eps \\
  \hline
 3 & 
 \bps 2 \eps & 
 \bps 1 & 3 \eps & 
 \bps 0 & 2 & 4 \eps & 
 \bps \cdots\cdots\cdots \eps & 
 \bps l-5 & l-3 & l-1 \eps & 
 \bps l-4 & l-2 \eps & 
 \bps l-3 \eps \\ 
 \hline
 \vdots & 
 \bps \vdots\\ \vdots\\ \vdots \eps & 
 \bps \vdots\\ \vdots\\ \vdots \eps & 
 \bps \vdots\\ \vdots\\ \vdots \eps & 
 \bps  && \\ & \ddots & \\ && \eps & 
 \bps \vdots\\ \vdots\\ \vdots \eps & 
 \bps \vdots\\ \vdots\\ \vdots \eps & 
 \bps \vdots\\ \vdots\\ \vdots \eps \\
  \hline 
 l-2 & 
 \bps l-3 \eps & 
 \bps l-4 & l-2 \eps & 
 \bps l-5 & l-3 & l-1 \eps & 
 \bps \cdots\cdots\cdots \eps & 
 \bps 0 & 2 & 4 \eps & 
 \bps 1 & 3 \eps & 
 \bps 2 \eps \\
  \hline
 l-1 & 
 \bps l-2 \eps & 
 \bps l-3 & l-1 \eps & 
 \bps l-4 & l-2 \eps & 
 \bps \cdots\cdots\cdots \eps & 
 \bps 1 & 3 \eps & 
 \bps 0 & 2 \eps & 
 \bps 1 \eps \\
  \hline
 l & 
 \bps l-1 \eps & 
 \bps l-2 \eps & 
 \bps l-3 \eps & 
 \bps \cdots\cdots\cdots \eps & 
 \bps 2 \eps & 
 \bps 1 \eps & 
 \bps 0 \eps \\
  \hline
 \end{array}
\end{equation*}

\vspace*{0.3cm}

\noindent
{\bf Type $D_l$}\ $(h-2=2l-4)$\ :
\begin{equation*}
 \renewcommand{\arraystretch}{2.5}
\begin{array}{@{\vrule width 1pt\ \ }c|c|c|c|c@{\ \vrule width 1pt}}
  \hline
 k\backslash k' & 1 & k'\ (\bps 2\le k'\le l-2\eps ) & l-1 & l \\ 
  \hline
 1 & 
 \bps 0 & 2l-4 \eps & 
 \bps k'-1 & 2l-3-k' \eps & 
 \bps l-2 \eps & 
 \bps l-2 \eps \\
  \hline
\underset{(2\le k\le l-2)}{k} & 
 \bps k-1\\ 2l-3-k \eps & 
 \bps |k'-k| \ \, |k'-k|+2 \ \, \cdots\cdots \quad  \\
 \quad \cdots\cdots \ \, k+k'-4 \ \, k+k'-2\ , \\
 2l-2-(k+k') \ 2l-(k+k') \ \cdots\cdots\quad \\ 
 \quad\cdots\cdots \ \, 2l-6-|k'-k| \ \, 2l-4-|k'-k| \\ \eps & 
 \bps l-1-k & l+1+k & \cdots \\ 
      \cdots & l-5+k & l-3+k \eps & 
 \bps l-1-k & l+1+k & \cdots \\ 
      \cdots & l-5+k & l-3+k \eps \\
  \hline
 l-1 & 
 \bps l-2 \eps & 
 \bps l-1-k & l+1+k & \cdots \\ 
      \cdots & l-5+k & l-3+k \eps & 
 \bps 0 & 4 & 8 & \cdots & 2l-4 \\ 
        &   &   &       (l:& \text{even}) \\ 
      0 & 4 & 8 & \cdots & 2l-6 \\ 
        &   &   &       (l:& \text{odd}) 
 \eps & 
 \bps 2 & 6 & 10 & \cdots & 2l-6 \\ 
         &   &   &       (l:& \text{even}) \\ 
      2 & 6 & 10 & \cdots & 2l-4 \\ 
        &   &   &       (l:& \text{odd}) 
      \eps \\ 
  \hline 
 l & 
 \bps l-2 \eps & 
 \bps l-1-k & l+1+k & \cdots \\ 
      \cdots & l-5+k & l-3+k \eps & 
 \bps 2 & 6 & 10 & \cdots & 2l-6 \\ 
         &   &   &       (l:& \text{even}) \\ 
      2 & 6 & 10 & \cdots & 2l-4 \\ 
        &   &   &       (l:& \text{odd}) 
 \eps & 
 \bps 0 & 4 & 8 & \cdots & 2l-4 \\
        &   &   &       (l:& \text{even}) \\ 
      0 & 4 & 8 & \cdots & 2l-6 \\ 
        &   &   &       (l:& \text{odd}) 
 \eps \\
  \hline
 \end{array}
\end{equation*}

\vspace*{0.3cm}

\noindent
{\bf Type $E_6$}\ $(h-2=10)$\ :
\begin{equation*}
 \renewcommand{\arraystretch}{1.5}
\begin{array}{@{\vrule width 1pt\ \ }c|c|c|c|c|c|c@{\ \vrule width 1pt}}
  \hline
 k\backslash k' & 1 & 2 & 3 & 4 & 5 & 6 \\ 
  \hline
 1 & 
 \bps 0 & 4 & 6 & 10 \eps & 
 \bps 1 & 3 & 5^2 & 7 & 9 \eps & 
 \bps 2 & 4 & 6 & 8 \eps & 
 \bps 2 & 4 & 6 & 8 \eps & 
 \bps 3 & 7 \eps & 
 \bps 3 & 7 \eps \\
  \hline
 2 & 
  \bps 1 & 3 & 5^2 & 7 & 9 \eps & 
 \bps 0 & 2^2 & 4^3 & 6^3 & 8^2 & 10 \eps & 
 \bps 1 & 3^2 & 5^2 & 7^2 & 9 \eps & 
 \bps 1 & 3^2 & 5^2 & 7^2 & 9 \eps & 
 \bps 2 & 4 & 6 & 8 \eps & 
 \bps 2 & 4 & 6 & 8 \eps \\
  \hline
 3 & 
  \bps 2 & 4 & 6 & 8 \eps & 
 \bps 1 & 3^2 & 5^2 & 7^2 & 9 \eps & 
 \bps 0 & 2 & 4 & 6^2 & 8 \eps & 
 \bps 2 & 4^2 & 6 & 8 & 10 \eps & 
 \bps 1 & 5 & 7 \eps & 
 \bps 3 & 5 & 9 \eps \\
  \hline 
 4 & 
 \bps 2 & 4 & 6 & 8 \eps & 
 \bps 1 & 3^2 & 5^2 & 7^2 & 9 \eps & 
 \bps 2 & 4^2 & 6 & 8 & 10 \eps & 
 \bps 0 & 2 & 4 & 6^2 & 8 \eps & 
 \bps 3 & 5 & 9 \eps & 
 \bps 1 & 5 & 7 \eps \\
  \hline
 5 & 
 \bps 3 & 7 \eps & 
 \bps 2 & 4 & 6 & 8 \eps & 
 \bps 1 & 5 & 7 \eps & 
 \bps 3 & 5 & 9 \eps & 
 \bps 0 & 6 \eps & 
 \bps 4 & 10 \eps \\
  \hline
 6 & 
 \bps 3 & 7 \eps & 
 \bps 2 & 4 & 6 & 8 \eps & 
 \bps 3 & 5 & 9 \eps & 
 \bps 1 & 5 & 7 \eps & 
 \bps 4 & 10 \eps & 
 \bps 0 & 6 \eps  \\
  \hline 
 \end{array}
\end{equation*}

\vspace*{0.3cm}

\noindent 
{\bf Type $E_7$}\ $(h-2=16)$\ :
\begin{equation*}
 \renewcommand{\arraystretch}{2.0}
\begin{array}{@{\vrule width 1pt\ \ }c|c|c|c|c|c|c|c@{\ \vrule width 1pt}}
  \hline
 k\backslash k' & 1 & 2 & 3 & 4 & 5 & 6 & 7 \\ 
  \hline
 1 & 
 \bps 0 & 6 & 10 & 16 \\ &&& \eps & 
 \bps 1 & 5 & 7 \\ 9 & 11 & 15\\ && \eps & 
 \bps 2 & 4 & 6 \\ &8^2& \\ 10 & 12 & 14\\ && \eps & 
 \bps 3 & 7 & 9 & 13\\ &&& \eps & 
 \bps 3 & 5 & 7 \\ 9 & 11 & 13\\ &&& \eps & 
 \bps 4 & 6 & 10 & 12\\ &&& \eps & 
 \bps 5 & 11\\ & \eps \\
  \hline
 2 & 
 \bps 1 & 5 & 7 \\ 9 & 11 & 15\\ && \eps & 
 \bps 0 & 2 & 4\\ 6^2 & 8^2 & 10^2\\ 12 & 14 & 16\\ && \eps & 
 \bps 1 & 3^2 & 5^2 & 7^3 \\ 9^3 & 11^2 & 13^2 & 15\\ &&& \eps & 
 \bps 2 & 4 & 6\\ &8^2& \\ 10 & 12 & 14\\ && \eps & 
 \bps 2 & 4^2 & 6^2\\ &8^2& \\ 10^2 & 12^2 & 14\\ && \eps & 
 \bps 3 & 5^2 & 7 \\ 9 & 11^2 & 13\\ && \eps & 
 \bps 4 & 6 & 10 & 12\\ &&& \eps \\
  \hline
 3 & 
 \bps 2 & 4 & 6\\ &8^2& \\ 10 & 12 & 14\\ && \eps & 
 \bps 1 & 3^2 & 5^2 & 7^3 \\ 9^3 & 11^2 & 13^2 & 15\\ &&& \eps & 
 \bps 0 & 2^2 & 4^3 \\ 6^4 & 8^4 & 10^4 \\ 12^3 & 14^2 & 16\\ && \eps & 
 \bps 1 & 3 & 5^2 & 7^2 \\ 9^2 & 11^2 & 13 & 15\\ &&& \eps & 
 \bps 1 & 3^2 & 5^3 & 7^3 \\ 9^3 & 11^3 & 13^2 & 15\\ &&& \eps & 
 \bps 2 & 4^2 & 6^2\\ & 8^2& \\ 10^2 & 12^2 & 14\\ && \eps & 
 \bps 3 & 5 & 7 \\ 9 & 11 & 13\\ && \eps \\
  \hline 
 4 & 
 \bps 3 & 7 & 9 & 13\\ &&& \eps & 
 \bps 2 & 4 & 6\\ &8^2& \\ 10 & 12 & 14\\ && \eps & 
 \bps 1 & 3 & 5^2 & 7^2 \\ 9^2 & 11^2 & 13 & 15\\ &&& \eps & 
 \bps 0 & 4 & 6 \\ &8& \\ 10 & 12 & 16\\ && \eps & 
 \bps 2 & 4 & 6^2\\ &8& \\ 10^2 & 12 & 14\\ && \eps & 
 \bps 3 & 5 & 7 \\ 9 & 11 & 13\\ && \eps & 
 \bps 4 & 8 & 12\\ && \eps \\
  \hline
 5 & 
 \bps 3 & 5 & 7 \\ 9 & 11 & 13\\ && \eps & 
 \bps 2 & 4^2 & 6^2\\ &8^2& \\ 10^2 & 12^2 & 14\\ && \eps & 
 \bps 1 & 3^2 & 5^3 & 7^3 \\ 9^3 & 11^3 & 13^2 & 15\\ &&& \eps & 
 \bps 2 & 4 & 6^2\\ &8& \\ 10^2 & 12 & 14\\ && \eps & 
 \bps 0 & 2 & 4^2 \\ 6^2 & 8^3 & 10^2 \\ 12^2 & 14 & 16\\ && \eps & 
 \bps 1 & 3 & 5 & 7^2 \\ 9^2 & 11 & 13 & 15\\ &&& \eps & 
 \bps \, 2\ \ 6\\ 8 \\ 10\ 14\\ \ \eps \\
  \hline
 6 & 
 \bps 4 & 6 & 10 & 12\\ &&& \eps & 
 \bps 3 & 5^2 & 7 \\ 9 & 11^2 & 13\\ && \eps & 
 \bps 2 & 4^2 & 6^2\\ & 8^2& \\ 10^2 & 12^2 & 14\\ && \eps & 
 \bps 3 & 5 & 7 \\ 9 & 11 & 13\\ && \eps & 
 \bps 1 & 3 & 5 & 7^2 \\ 9^2 & 11 & 13 & 15\\ &&& \eps & 
 \bps 0 & 2 & 6\\ &8^2& \\ 10 & 14 & 16\\ && \eps & 
 \bps 1 & 7 & 9 & 15\\ &&& \eps \\
  \hline 
 7 & 
 \bps 5 & 11\\ & \eps & 
 \bps 4 & 6 & 10 & 12\\ &&& \eps & 
 \bps 3 & 5 & 7 \\ 9 & 11 & 13\\ && \eps & 
 \bps 4 & 8 & 12\\ && \eps & 
 \bps \, 2\ \ 6\\ 8 \\ 10\ 14\\ \ \eps & 
 \bps 1 & 7 & 9 & 15\\ &&& \eps & 
 \bps 0 & 8 & 16\\ && \eps \\
  \hline
\end{array}
\end{equation*}

\vspace*{0.3cm}

\noindent
{\bf Type $E_8$}\ $(h-2=28)$\ :
\begin{equation*}
 \renewcommand{\arraystretch}{2.5}
\begin{array}{@{\vrule width 1pt\ \ }c|c|c|c|c|c|c|c|c@{\ \vrule width 0.1pt}}
  \hline
 k\backslash k' & 1 & 2 & 3 & 4 & 5 & 6 & 7 & 8\\ 
  \hline
 1 & 
 \bps 0 & 10 \\ 18 & 28 \\ & \eps & 
 \bps 1 & 9 & 11 \\ 17 & 19 & 27 \\ & & \eps & 
 \bps \ 2\ \ 8 \ \, 10 \\ 12\ 16 \\ 18\ 20\ 26 \\ \\ \eps & 
 \bps \, 3 \ \ 7 \ \ 9 \\ 11\ 13 \\ 15\ 17 \\ 19\ 21\ 25  \\ \\ \eps & 
 \bps \, 4\ \ 6\\ \, 8\ \, 10\\ 
 12\ 14^2\, 16\\ 18\ 20\\ 22\ 24  \\ \\ \eps & 
 \bps 5 & 9 & 13 \\ 15 & 19 & 23 \\ \\ \eps & 
 \bps \, 5\ \ 7 \ \ 11 \\ 13\ 15 \\ 17\ 21\ 23  \\ \\ \eps & 
 \bps 6 & 12 \\ 16 & 22 \\ \\ \eps \\
  \hline
 2 & 
 \bps 1 & 9 & 11 \\ 17 & 19 & 27 \\ \eps & 
\bps 0\ \ 2\ \ 8 \\ 10^2\ \, 12\\ 
 \, 16\ \ \, 18^2 \\ 20\ 26\ 28 \\ \\ \eps & 
 \bps 1\! & 3\! & 7 \\ 9^2\! & 11^2\! & 13 \\ 
 15\! & 17^2\! & 19^2 \\ 21\! & 25\! & 27 \\ \\ \eps & 
 \bps  2\ \ 4\ \ 6 \\ 8^2\ \, 10^2 \\ 12^2\, 14^2\, 16^2 \\ 
       18^2\ 20^2 \\ 22\ 24\ 26 \\ \\ \eps & 
 \bps 3\! & 5^2\! & 7^2\\ 9^2\! & 11^2\! & 13^3 \\ 
      15^3\! & 17^2\! & 19^2\! \\ 21^2\! & 23^2\! & 25 \\ \\ \eps & 
 \bps \, 4\ \ 6 \\ \, 8\ \, 10 \\ 12\ 14^2\, 16 \\ 
 18\ 20 \\ 22\ 24 \\ \\ \eps & 
 \bps \ 4\ \ 6^2 \\ \, 8\ \, 10 \\ 12^2\, 14^2\, 16^2 \\ 
 18\ 20 \\ 22^2\ 24 \\ \\ \eps &
 \bps \, 5\ \ 7\ \ 11 \\ 13\ 15 \\ 17\ 21\ 23 \\ \\ \eps \\
  \hline
 3 & 
 \bps \ 2\ \ 8\ \, 10 \\ 12\ 16 \\ 18\ 20\ 26 \\ \\ \eps & 
 \bps 1\! & 3\! & 7 \\ 9^2\! & 11^2\! & 13 \\ 
 15\! & 17^2\! & 19^2 \\ 21\! & 25\! & 27 \\ \\ \eps & 
 \bps 0\! & 2\! & 4 \\ 6\! & 8^2\! & 10^3 \\ 12^2\! & 14^2\! & 16^2 \\
      18^3\! & 20^2\! & 22 \\ 24\! & 26\! & 28 \\ \\ \eps & 
 \bps \, 1\ \ \, 3\ \ 5^2 \\ 7^2\ \, 9^2\ 11^3 \\ 13^3\ 15^3 \\ 
      17^3\, 19^3\, 21^2 \\ 23^2\, 25\ 27\ \\ \\ \eps & 
 \bps 2\! & 4^2\! & 6^3 \\ 8^3\! & 10^3\! & 12^4 \\ 
      16^4\! & 18^3\! & 20^3 \\ 22^3\! & 24^2\! & 26 \\ \\ \eps & 
 \bps 3\! & 5\! & 7^2 \\ 9\! & 11^2\! & 13^2 \\ 15^2\! & 17^2\! & 19 \\
      21^2\! & 23\! & 25 \\ \\ \eps & 
 \bps 3\! & 5^2\! & 7^2 \\ 9^2\! & 11^2\! & 13^3 \\ 
      15^3\! & 17^2\! & 19^2 \\ 21^2\! & 23^2\! & 25 \\ \\ \eps & 
 \bps \, 4\ \ 6 \\ \ 8\ \, 10 \\ 12\ 14^2\, 16 \\ 
 18\ 20 \\ 22\ 24 \\ \\ \eps \\
  \hline 
 4 & 
 \bps \,3\ \ 7\ \ 9 \\ 11\ 13 \\ 15\ 17 \\ 19\ 21\ 25 \\ \\ \eps & 
 \bps \, 2\ \  4\ \ 6 \\ 8^2\ 10^2 \\ 12^2\, 14^2\, 16^2 \\ 
       18^2\ 20^2 \\ 22\ 24\ 26 \\ \\ \eps & 
 \bps \, 1\ \ \, 3\ \ 5^2 \\ 7^2\ \, 9^2\ 11^3 \\ 13^3\ 15^3 \\ 
      17^3\, 19^3\, 21^2 \\ 23^2\, 25\ 27\ \\ \\ \eps & 
 \bps 0\! & 2\! & 4^2 \\ 6^3\! & 8^3\! & 10^4 \\ 
      12^4\! & 14^4\! & 16^4 \\ 
      18^4\! & 20^3\! & 22^3 \\ 24^2\! & 26\! & 28 \\ \\ \eps & 
 \bps \, 1\ \ 3^2\ \, 5^3 \\ \, 7^4\ \, 9^4\ 11^5 \\ 13^5\ 15^5 \\ 
      17^5\, 19^4\, 21^4 \\ 23^3\, 25^2\, 27\ \\ \\ \eps & 
 \bps \ 2\ \  4\ \ 6^2 \\ \, 8^2\ 10^2 \\
      12^3\, 14^2\, 16^3 \\ 
      18^2\ 20^2 \\ 22^2\, 24\ 26 \\ \\ \eps & 
 \bps \ 2\ \ \, 4^2\ \, 6^2 \\ \, 8^3\ \, 10^3 \\ 
      12^3\, 14^4\, 16^3 \\ 
      18^3\, 20^3\\ 22^2\, 24^2\, 26\, \\ \\ \eps & 
 \bps 3\! & 5\! & 7 \\ 9^2\! & 11\! & 13^2 \\ 
      15^2\! & 17\! & 19^2 \\ 21\! & 23\! & 25 \\ \\ \eps \\
  \hline
 5 & 
 \bps \, 4\ \ 6\\ \, 8\ \, 10\\ 
 12\ 14^2\, 16\\ 18\ 20\\ 22\ 24  \\ \\ \eps & 
 \bps 3\! & 5^2\! & 7^2 \\ 9^2\! & 11^2\! & 13^3 \\ 
      15^3\! & 17^2\! & 19^2 \\ 21^2\! & 23^2\! & 25 \\ \\ \eps & 
 \bps 2\! & 4^2\! & 6^3 \\ 8^3\! & 10^3\! & 12^4 \\ 
      16^4\! & 18^3\! & 20^3 \\ 22^3\! & 24^2\! & 26 \\ \\ \eps & 
 \bps \, 1\ \ 3^2\ \, 5^3 \\ \, 7^4\ \, 9^4\ 11^5 \\ 13^5\ 15^5 \\ 
      17^5\, 19^4\, 21^4 \\ 23^3\, 25^2\, 27\ \\ \\ \eps & 
 \bps 0\! & 2^2\! & 4^3 \\ 6^4\! & 8^5\! & 10^6 \\ 12^6\! & 14^6\! & 16^6 \\
      18^6\! & 20^5\! & 22^4 \\ 24^3\! & 26^2\! & 28 \\ \\ \eps & 
 \bps \, 1\ \ 3\ \ 5^2 \\ \, 7^2\ 9^3\ 11^3 \\ 13^3\ 15^3 \\
      17^3\, 19^3\, 21^2 \\ 23^2\, 25\ 27\ \\ \\ \eps & 
 \bps \, 1\ \, 3^2\ \, 5^2 \\ \, 7^3\ 9^4\ 11^4 \\
      13^4\ 15^4 \\ 
      17^4\, 19^4\, 21^3 \\ 23^2\, 25^2\, 27\ \\ \\ \eps & 
 \bps \ 2\ \ 4\ \ 6 \\ \, 8^2\ 10^2 \\ 
      12^2\, 14^2\, 16^2 \\ 
      18^2\ 20^2 \\ 22\ 24\ 26 \\ \\ \eps \\
  \hline
 6 & 
 \bps 5 & 9 & 13 \\ 15 & 19 & 23 \\ \\ \eps & 
 \bps \, 4\ \ 6 \\ \ 8\ \, 10 \\ 12\ 14^2\, 16\\ 
  18\ 20 \\ 22\ 24 \\ \\ \eps & 
 \bps 3\! & 5\! & 7^2 \\ 9\! & 11^2\! & 13^2 \\ 15^2\! & 17^2\! & 19 \\
      21^2\! & 23\! & 25 \\ \\ \eps & 
 \bps \ 2\ \ 4\ \ 6^2 \\ 8^2\ 10^2 \\
      12^3\, 14^2\, 16^3 \\ 
      18^2\ 20^2 \\ 22^2\, 24\ 26 \\ \\ \eps & 
 \bps \, 1\ \ 3\ \ 5^2 \\ \, 7^2\ 9^3\ 11^3 \\ 13^3\ 15^3 \\
      17^3\, 19^3\, 21^2 \\ 23^2\, 25\ 27\ \\ \\ \eps & 
 \bps \, 0\ \ 4\ \ 6 \\ \ 8\ \ 10^2 \\ 12\ 14^2\, 16 \\ 
      18^2\ 20 \\ 22\ 24\ 28 \\ \\ \eps & 
 \bps \, 2\ \ 4\ \ 6 \\ 8^2\ 10^2 \\ 12^2\, 14^2\, 16^2 \\ 
      18^2\ 20^2 \\ 22\ 24\ 26 \\ \\ \eps & 
 \bps 3\ \ 7 \\ \, 9\ \, 11\ \, 13 \\ 15\ 17\ 19 \\ 21\  25 \\ \\ \eps \\
  \hline 
 7 & 
 \bps \, 5\ \ 7\ \ 11 \\ 13\ 15 \\ 17\ 21\ 23 \\ \\ \eps & 
 \bps \ 4\ \ \ 6^2 \\  8\ \ 10 \\ 12^2\, 14^2\, 16^2 \\ 
      18\ 20 \\ 22^2\ 24 \\ \\ \eps &
 \bps 3\! & 5^2\! & 7^2 \\ 9^2\! & 11^2\! & 13^3 \\ 
      15^3\! & 17^2\! & 19^2 \\ 21^2\! & 23^2\! & 25 \\ \\ \eps & 
 \bps \, 2\ \ \, 4^2\ \, 6^2 \\ \, 8^3\ \, 10^3 \\ 
      12^3\, 14^4\, 16^3 \\ 
      18^3\ 20^3 \\ 22^2\, 24^2\, 26 \\ \\ \eps & 
 \bps \, 1\ \, 3^2\ \, 5^2 \\ \, 7^3\ 9^4\ 11^4 \\
      13^4\ 15^4 \\ 
      17^4\, 19^4\, 21^3 \\ 23^2\, 25^2\, 27\ \\ \\ \eps & 
 \bps \, 2\ \ 4\ \ 6 \\ \, 8^2\ 10^2 \\ 12^2\, 14^2\, 16^2 \\ 
      18^2\ 20^2 \\ 22\ 24\ 26 \\ \\ \eps & 
 \bps 0\! & 2\! & 4 \\ 6^2\! & 8^2\! & 10^3 \\ 
      12^3\! & 14^2\! & 16^3 \\ 18^3\! & 20^2\! & 22^2 \\ 
      24\! & 26\! & 28 \\ \\ \eps & 
 \bps 1\! & 5\! & 7 \\ 9\! & 11^2\! & 13 \\ 15\! & 17^2\! & 19 \\ 
 21\! & 23\! & 27 \\ \\ \eps \\
  \hline
 8 & 
 \bps 6 & 12 \\ 16 & 22 \\ \\ \eps & 
 \bps \, 5\ \ 7\ \ 11 \\ 13\ 15 \\ 17\ 21\ 23 \\ \\ \eps & 
 \bps \, 4\ \ \, 6 \\ \ 8\ \, 10 \\ 12\  14^2\, 16 \\ 18\ 20 \\ 22\ 24 
      \\ \\ \eps & 
 \bps 3\! & 5\! & 7 \\ 9^2\! & 11\! & 13^2 \\ 
      15^2\! & 17\! & 19^2 \\ 21\! & 23\! & 25 \\ \\ \eps & 
 \bps \, 2 \ \ 4 \ \ 6 \\ 8^2 \ \, 10^2 \\ 
      12^2\, 14^2\, 16^2 \\ 
      18^2 \ \, 20^2 \\ 22 \ 24 \ 26 \\ \\ \eps & 
 \bps 3\ \ 7 \\ \, 9\ \, 11\ \, 13 \\ 15\ 17\ 19 \\ 21\  25 \\ \\ \eps &
 \bps 1\! & 5\! & 7 \\ 9\! & 11^2\! & 13 \\ 15\! & 17^2\! & 19 \\ 
 21\! & 23\! & 27 \\ \\ \eps & 
 \bps \, 0\ \ 6\ \, 10 \\ 12\ 16 \\ 18\ 22\ 28 \\ \\ \eps \\
  \hline
\end{array}
\end{equation*}

\vspace*{0.5cm}

\noindent
For each ADE case, 
one can easily check the followings. 
\begin{enumerate}
 \item 
One has $\fC(k,k')=\fC(k',k)$ for any $k,k'\in\Pi$. 
This implies, for any $k,k'\in\Pi$, 
$\HomDz{R}(\ti{M}^k_n,\ti{M}^{k'}_{n'})\simeq
\HomDz{R}(\ti{M}^{k'}_{n'},\ti{M}^k_{n''})$ holds for some 
$n,n,n''\in\Z$ 
such that $\phi^{k}_{n''}-\phi^{k'}_{n'}=\phi^{k'}_{n'}-\phi^k_n$. 

 \item A consequence of the Serre duality 
(Theorem \ref{thm:grading} (iii)), 
\begin{equation*}
\dim_\C(\HomDz{R}(\ti{M}^k_{n},\ti{M}^{k'}_{n'}))
=\dim_\C(\HomDz{R}(\ti{M}^{k'}_{n'},\S(\ti{M}^k_{n})))\ ,
\end{equation*} 
can be checked as follows. 
For an indecomposable object $\ti{M}^k_n\in\tHMF{R}$, 
$k^\S\in\Pi$ denotes the vertex such that 
$[F(\S(\ti{M}^k_n))]=[M^{k^\S}]\in\Pi$. 
Then, for the given multi-set $\fC(k,k')$, one has 
$\fC(k',k^\S)=\{h-2-c\ |\ c\in\fC(k,k')\}$.

 \item One can check Theorem \ref{thm:grading} (ii-b): 
for any two indecomposable objects 
$\ti{M}^k_{n},\ti{M}^{k'}_{n'}\in\tHMF{R}$, 
$\dim_\C(\HomDz{R}(\ti{M}^k_{n},\ti{M}^{k'}_{n'}))=1$ 
if $h|\phi^{k'}_{n'}-\phi^k_{n}|=d(k,k')$. 
Correspondingly, one has 
$\sharp\{c\in\fC(k,k')\ |\ c=d(\ti{M}^k_{n},\ti{M}^{k'}_{n'})\}=1$. 
In addition, the Serre duality implies that 
$\dim_\C(\HomDz{R}(\ti{M}^k_{n},\ti{M}^{k'}_{n'}))=1$ 
if $h|\phi^{k'}_{n'}-\phi^k_{n}|=h-2-d(k^\S,k')$. 
Namely, one has 
$\sharp\{c\in\fC(k,k')\ |\
c=h-2-d(k^\S,k')\}=1$. 
These facts, together with Corollary \ref{cor:serre1}, imply that 
$\fC(k,k')=\fC(k',k)$ is described in the following form: 
\begin{equation*}
 \qquad \fC(k,k')=\{c_1,\cdots, c_s\}\ ,\qquad 
 d(k,k')=c_1<c_2\le\cdots\le c_{s-1}< c_s=h-2-d(k^\S,k)\ 
\end{equation*}
for some $s\in\Z_{>0}$. 
\end{enumerate}
\end{tbl}

 \appendix

\section{Another proof of Theorem \ref{thm:main} \ by Kazushi Ueda}

In this appendix, we give another proof of Theorem \ref{thm:main} 
which avoids the use of the classification of
Cohen--Macaulay modules
on simple singularities
due to Auslander 
%\cite{Auslander_RSASS},
\cite{a:1}
and is based on the theory of weighted projective lines
by Geigle and Lenzing \cite{Geigle-Lenzing_WPC, Geigle-Lenzing_PC}
and a theorem of Orlov 
%\cite{Orlov_DCCSTCS}
\cite{o:2} 
on triangulated categories of graded B-branes.

\subsection{Weighted projective lines of Geigle and Lenzing}

Let $k$ be a field.
For a sequence $\pp = (p_0, p_1, p_2)$
of nonzero natural numbers,
let $L(\pp)$
be the abelian group of rank one
generated by four elements
$\vx_0, \vx_1, \vx_2, \vc$
with relations
$p_0 \vx_0 = p_1 \vx_1 = p_2 \vx_2 = \vc$,
and consider the $k$-algebra
\begin{equation*}
 R(\pp):=k[x_0,x_1,x_2] / (x_0^{p_0} - x_1^{p_1} + x_2^{p_2})
\end{equation*}
graded by
$\deg(x_s)=\vx_s \in L(\pp)$ for $s=0,1,2$. 
Define the category
%$\qgr R(\pp)$
of coherent sheaves on the weighted projective line
of weight $\pp$
as the quotient category
\begin{equation*}
 \qgr R(\pp):=\gr R(\pp) / \tor R(\pp)
\end{equation*}
of the abelian category $\gr R(\pp)$
of finitely-generated
$L(\pp)$-graded $R(\pp)$-modules
by its full subcategory $\tor R(\pp)$ 
consisting of torsion modules. 
This definition is equivalent to the original one by Geigle and 
Lenzing due to Serre's theorem
in \cite[section 1.8]{Geigle-Lenzing_WPC}.
Let $\pi : \gr R(\pp) \to \qgr R(\pp)$ be the natural projection.
For $M \in \gr R(\pp) $ and $\vx\in L(\pp)$, 
let $M(\vx)$ be the graded $R(\pp)$-module 
obtained by shifting the grading by $\vx$;
$
 M(\vx)_{\vec{n}}=M_{\vec{n}+\vec{x}},
$
and put
$
 \scO(\vec{n}) = \pi R(\pp)(\vec{n}).
$
Then
%for any weight sequence $\pp$,
the sequence
%\begin{eqnarray*}
% & (E_0, \ldots, E_N)
%  = \left( \pi R(\pp),
%     \pi R(\pp)(\vecx_0),
%     \pi R(\pp)(2 \vecx_0), \cdots,
%     \pi R(\pp)((p_0 - 1) \vecx_0),
%    \right. \\ & \qquad \left.
%     \pi R(\pp)(\vecx_1),
%     \pi R(\pp)(2 \vecx_1), \cdots,
%     \pi R(\pp)((p_{2} - 1) \vecx_2),
%     \pi R(\pp)(\vecc) \right)
%\end{eqnarray*}
\begin{eqnarray*}
 & (E_0, \ldots, E_N)
  = \left( \scO,
     \scO(\vecx_0),
     \scO(2 \vecx_0), \cdots,
     \scO((p_0 - 1) \vecx_0),
    \right. \\ & \qquad \left.
     \scO(\vecx_1),
     \scO(2 \vecx_1), \cdots,
     \scO((p_{2} - 1) \vecx_2),
     \scO(\vecc) \right)
\end{eqnarray*}
of objects of $\qgr R(\pp)$,
where $N = p_0 + p_1 + p_2 - 2$,
is a full strong exceptional collection
by \cite[Proposition 4.1]{Geigle-Lenzing_WPC}.
Define the {\em dualizing element}
$\vecomega \in L(\pp)$ by
$$
 \vecomega = \vecc - \vecx_0 - \vecx_1 - \vecx_2
$$
and a $\bZ$-graded subalgebra $R'(\pp)$ of $R(\pp)$ by
\begin{equation} \label{eq:Rprime}
 R'(\pp) = \bigoplus_{n \geq 0} R'(\pp)_n, \qquad
 R'(\pp)_n = R(\pp)_{-n \vecomega}.
\end{equation}
A weight sequence $\pp = (p_0, p_1, p_2)$
is called {\em of Dynkin type}
if
$$
 \frac{1}{p_0} + \frac{1}{p_1} + \frac{1}{p_2} > 1.
$$
A weight sequence of Dynkin type
yields
the rational double point of the corresponding type:

\begin{prop}[{\cite[Proposition 8.4.]{Geigle-Lenzing_PC}}]
\label{prop:Rprime}
For a weight sequence $\pp$ of Dynkin type,
the $\bZ$-graded algebra $R'(\pp)$ has the form
$$
 k[x, y, z] / f_{\pp}(x, y, z)
$$
where 
$f_{\pp}(x, y, z)$ is the polynomial of type $A_{p+q}$ if 
$\pp=(1, p, q)$, $D_{2 l - 2}$ if $\pp=(2, 2, 2 l)$, 
$D_{2 l - 1}$ if $\pp=(2, 2, 2 l + 1)$, 
$E_6$ if $\pp=(2, 3, 3)$, 
$E_7$ if $\pp=(2, 3, 4)$, and 
$E_8$ if $\pp=(2, 3, 5)$. 
%
%the homogeneous generators $(x, y, z)$ 
%and the relation $f_{\pp}(x, y, z)$ are displayed
%in the following table:
%$$
%\begin{array}[t]{c|cccc}
% &
% \text{weight} &
% \text{Generators $(x, y, z)$} &
% \text{$\bZ$-degrees} & \text{Relations $f_{\pp}$} \\
% \hline
% A_{p+q} & (1, p, q) & (x_1 x_2, x_1^{p+q}, x_2^{p+q}) &
% (1, p, q) & x^{p+q} - y z \\
% D_{2 l - 2} & (2, 2, 2 l) & (x_2^2, x_0^2, x_0 x_1 x_2) &
% (2, 2 l, 2 l + 1) & z^2 + x (y^2 + y x^l) \\
% D_{2 l - 1} & (2, 2, 2 l + 1) & (x_2^2, x_0 x_1, x_0^2 x_2) &
% (2, 2 l + 1, 2 l + 2) & z^2 + x (y^2 + z x^l) \\
% E_6 & (2, 3, 3) & (x_0, x_1 x_2, x_1^3) &
% (3, 4, 6) & z^2 + y^3 + x^2 z \\
% E_7 & (2, 3, 4) & (x_1, x_2^2, x_0 x_2) &
% (4, 6, 9) & z^2 + y^3 + x^3 y \\
% E_8 & (2, 3, 5) & (x_2, x_1, x_0) &
% (6, 10, 15) & z^2 + y^3 + x^5
%\end{array}
%$$
\end{prop}
Moreover, we have the following:
\begin{prop}[{\cite[Proposition 8.5]{Geigle-Lenzing_PC}}]
\label{prop:cohX}
For a weight sequence $\pp$ of Dynkin type,
there exists a natural equivalence
$$
 \qgr R(\pp) \cong \qgr R'(\pp),
$$
where $\qgr R'(\pp)$
is the quotient category of
the abelian category $\gr R'(\pp)$
of finitely-generated
$\Z$-graded $R'(\pp)$-modules
by its full subcategory $\tor R'(\pp)$ 
consisting of torsion modules. 
\end{prop}

 \subsection{Graded B-branes on simple singularities}

For the $\bZ$-graded algebra $R'(\pp)$
given above,
define the {\em triangulated category of graded B-branes}
as the quotient category
$$
 \Dbsing(R'(\pp)) := D^b(\gr R'(\pp)) / D^b(\grproj R'(\pp))
$$
of the bounded derived category $D^b(\gr R'(\pp))$
by its full triangulated subcategory $D^b(\grproj R'(\pp))$
consisting of perfect complexes,
i.e., bounded complexes of projective $\bZ$-graded modules
%\cite{Orlov_DCCSTCS}.
\cite{o:2}.
$\Dbsing(R'(\pp))$ is equivalent to
$HMF^{gr}_{k[x,y,z]}(f_{\pp})$. 
%in the notation of \cite{Kajiura-Saito-Takahashi_MF}.
%Since 
Note $R'(\pp)$ is Gorenstein with Gorenstein parameter one,
i.e.,
\begin{equation*}
 \mathrm{Ext}^i_A(k,A)=
 \begin{cases}
   k(1)\ & \text{if}\ i=2,\\
    0  \ &\text{otherwise},\ 
 \end{cases} 
\end{equation*}
which follows from, e.g.,
\cite[Corollary 2.2.8 and Proposition 2.2.10]{Goto-Watanabe_OGR_I}.
Thus
there exists a fully faithful functor
$\Phi_0 : \Dbsing(R'(\pp)) \to D^b(\qgr R'(\pp))$
and a semiorthogonal decomposition
$$
 D^b(\qgr R'(\pp))
  = \langle E_0, \Phi_0 \Dbsing(R'(\pp)) \rangle
$$
by 
%\cite[Theorem 2.5.(i)]{Orlov_DCCSTCS}.
\cite[Theorem 2.5.(i)]{o:2}.
Therefore,
$\Dbsing(R'(\pp))$ is equivalent
to the full triangulated subcategory of $D^b(\qgr R'(\pp))$
generated by the strong exceptional collection
$(E_1, \ldots, E_N)$.
Its total morphism algebra
$\End(\bigoplus_{i=1}^n E_i)$
is isomorphic to the path algebra
of the Dynkin quiver

\vspace*{0.45cm}

%\begin{equation*} %\label{eq:Dynkin_quiver}
%\vecDelta(\pp) :
%\begin{array}{c@{\hskip13mm}c@{\hskip13mm}
%c@{\hskip13mm}c@{\hskip13mm}c}
% \Rnode{01}{\vecx_0} &
% \Rnode{02}{2 \vecx_0} &
% \Rnode{03}{\cdots} &
% \Rnode{04}{(p_0 - 1)}\Rnode{05}{\vecx_0} & \\[5mm]
% \Rnode{11}{\vecx_1} &
% \Rnode{12}{2 \vecx_1} &
% \Rnode{13}{\cdots} &
% \Rnode{14}{(p_1 - 1 )\vecx_1} &
% \Rnode{c}{\vecc} \\[5mm]
% \Rnode{21}{\vecx_2} &
% \Rnode{22}{2 \vecx_2} &
% \Rnode{23}{\cdots} &
% \Rnode{24}{(p_2 - 1)} \Rnode{25}{\vecx_2} & 
%\end{array}
%\psset{nodesep=3pt}
%%\everypsbox{\scriptstyle}
%\ncline{->}{01}{02}\Aput{x_0}
%\ncline{->}{02}{03}\Aput{x_0}
%\ncline{->}{03}{04}\Aput{x_0}
%\ncline{->}{05}{c}\Aput{x_0}
%\ncline{->}{11}{12}\Aput{x_1}
%\ncline{->}{12}{13}\Aput{x_1}
%\ncline{->}{13}{14}\Aput{x_1}
%\ncline{->}{14}{c}\Aput{x_1}
%\ncline{->}{21}{22}\Aput{x_2}
%\ncline{->}{22}{23}\Aput{x_2}
%\ncline{->}{23}{24}\Aput{x_2}
%\ncline{->}{25}{c}\Bput{x_2}
%\end{equation*}
\begin{center}
 \includegraphics{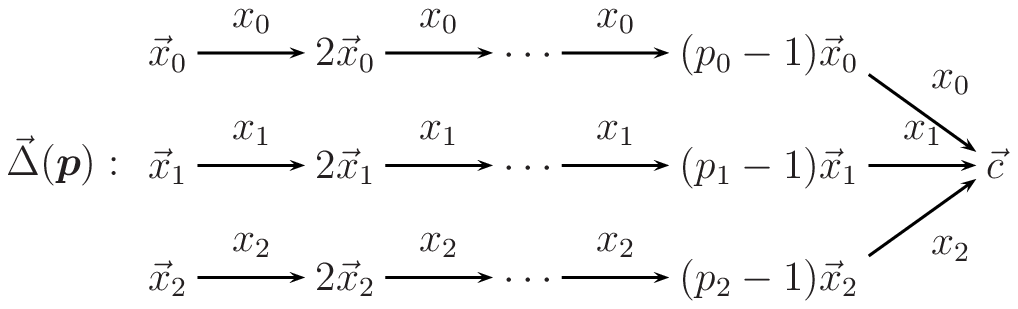}
\end{center}
of the corresponding type,
obtained from the quiver
appearing in \cite[section 4]{Geigle-Lenzing_WPC}
by removing the leftmost vertex.
Since
$\Dbsing(R'(\pp))$ is an enhanced triangulated category, 
the equivalence 
$
 \Dbsing(R'(\pp)) \simeq D^b(\mathrm{mod\text{--}}\C\vec{\Delta}(\pp))
$
follows from Bondal and Kapranov 
%\cite[Theorem 1]{Bondal-Kapranov_ETC}.
\cite[Theorem 1]{bk:1}.

%%%%%%%%%%%%%%%%%%%%%%%%%%%%%%%%%%%%%%%%%%%%%%%%%%%%%%%%%%%%%%%%%%%%%%%%%%%%%%

%
%%%%%%%%%%%%%%%%%%%%%%%%%%%%%%%%%%%%%%%%%%%%%%%%%%%%%%%%%%%%%%%%%%%%%%%%%%%%%%
%

\begin{thebibliography}{[AR]}
%
\bibitem[A]{a:1}
	M.~Auslander, 
	{\it Rational singularities and almost split sequences}, 
 	Trans. Amer. Math. Soc. 293 (1986) %no. 2, 
	511--531. 
\bibitem[AR1]{ar:ass}
	M.~Auslander, I.~Reiten, 
	{\it Representation theory of Artin algebras III,IV,V,VI}, 
	Comm. Algebra 3 (1975), 239--294; 5 (1977) 443--518; 
	5 (1977) 519--554; 6 (1978) 257--300. 
\bibitem[AR2]{ar:1}
	M.~Auslander and I.~Reiten,
	{\it Almost split sequences for rational double points},
	Trans. AMS., 302 (1987) 87-99.
\bibitem[AR3]{ar:2}
	M.~Auslander and I.~Reiten,
        {\it Cohen-Macaulay modules for graded Cohen-Macaulay rings 
   	and their completions}, 
	Commutative algebra (Berkeley, CA, 1987), 21--31, 
	Math. Sci. Res. Inst. Publ., 15, Springer, New York, 1989. 
\bibitem[Bd]{bd:1}
	T.~Bridgeland,
	{\it Stability conditions on triangulated category},
	math.AG/0212237.
\bibitem[Bs]{bs:1}
	E.~Brieskorn, 
	{\it Singular elements of semi-simple algebraic groups}, 
	Actes du Congr\`es International des Math\'ematiciens 
	(Nice, 1970), Tome 2, pp. 279--284. 
	Gauthier-Villars, Paris, 1971. 
\bibitem[BK]{bk:1}
	A.~Bondal and M.~Kapranov,
	{\it Enhanced triangulated categories},
	Math. USSR Sbornik, Vol.70, (1991) No.1, 93--107.
\bibitem[E]{e:1}
	D.~Eisenbud,
	{\it Homological algebra on a complete intersection, 
	with an application to group representations}, 
	Trans. AMS., 260 (1980) 35--64.
\bibitem[Ga]{g:1}
	P.~Gabriel,
	{\it Unzerleghare Derstellungen I},
	Manuscripta Math., {\bf 6} (1972) 71--163.
%
\bibitem[GL1]{Geigle-Lenzing_WPC}
       W.~Geigle and H.~Lenzing,
       {\it A class of weighted projective curves arising in
       representation theory of finite-dimensional algebras},
       Singularities, representation of algebras, and vector bundles 
       (Lambrecht, 1985), 9--34, 
        Lecture Notes in Math., 1273, Springer, Berlin, 1987. 
\bibitem[GL2]{Geigle-Lenzing_PC}
       W.~Geigle and H.~Lenzing,
       {\it Perpendicular categories with applications to
       representations and sheaves},
       J. Algebra, {\bf 144} %(2):273--343, 1991.
       (1991) 273--343. 
%
\bibitem[GM]{gm:1}
	I.~Gelfand and Y.~Manin,
	{\it Methods of Homological Algebra},
	Springer-Verlag (1994).
%
\bibitem[GW]{Goto-Watanabe_OGR_I}
       S.~Goto and K.~Watanabe,
       {\it On graded rings.{I}}, 
       J. Math. Soc. Japan, {\bf 30}%(2) 
       (1978) 179--213.
%
\bibitem[Gr]{gr}
	A.~Grothendieck, 
	{\it Groupes de classes des cat\'egories ab\'eliennes et 
	triangul\'es}' 
	Complexes parfaits, in: SGA 5, Lecture Notes in Math. 589, 
	Springer, 1977, 351--371. 
\bibitem[Ha]{h}
        D.~Happel, 
        {\it Triangulated categories in the representation theory of
        finite-dimensional algebras}, 
        London Mathematical Society Lecture Note Series, 119. 
        Cambridge University Press, Cambridge, 1988. x+208 pp. 
\bibitem[HW]{hw:1}
	K.~Hori and J.~Walcher,
	{\it F-term equations near Gepner points},
        JHEP {\bf 0501} (2005) 008, 
	hep-th/0404196. 
\bibitem[K]{k:1}
	H.~Kn\"{o}rrer,
	{\it Cohen-Macauley modules on hypersurface singularities I},
	Invent. Math., {\bf 88} (1987) 153-164.
\bibitem[KL1]{kl:1}
	A.~Kapustin and Y.~Li,
	{\it D-branes in Landau-Ginzburg models and algebraic geometry},
	JHEP {\bf 0312}, 005 (2003), hep-th/0210296.
\bibitem[KL2]{kl:2}
	A.~Kapustin and Y.~Li,
	{\it Topological Correlators in Landau-Ginzburg Models with Boundaries},
	Adv. Theoret. Math. Phys. {\bf 7} (2003) 727-749, hep-th/0305136.
\bibitem[KL3]{kl:3}
        A.~Kapustin and Y.~Li,
        {\it D-branes in topological minimal models: 
	The Landau-Ginzburg approach},
        JHEP {\bf 0407} (2004) 045, hep-th/0306001. 
\bibitem[KS]{ks:1}
	M.~Kashiwara and P.~Shapira,
	{\it Sheaves on Manifolds},
	Grundleheren 292, Springer-Verlag (1990).
\bibitem[O1]{o:1}
	D.~Orlov,
	{\it Triangulated Categories of Singularities and 
	D-branes in Landau-Ginzburg Models},
        Tr. Mat. Inst. Steklova 246 (2004), Algebr. Geom. Metody, 
  	Svyazi i Prilozh., 240--262; 
	translation in Proc. Steklov Inst. Math. 2004, no. 3 (246),
	227--248, 
	math.AG/0302304. 
\bibitem[O2]{o:2}
        D.~Orlov, 
        {\it Derived categories of coherent sheaves and triangulated
        categories of singularities}, 
        math.AG/0503632. 
\bibitem[R]{r:1}
        C.~M.~Ringel, 
	{\it Tame algebras and integral quadratic forms}, 
	Lecture Notes in Mathematics, {\bf 1099}. Springer-Verlag, Berlin,
	1984. xiii+376 pp. 
\bibitem[Sa1]{sa:weight}
	K.~Saito, 
	{\it Regular system of weights and associated singularities}, 
 	Complex analytic singularities, 479--526, 
	Adv. Stud. Pure Math., 8, North-Holland, Amsterdam, 1987. 
\bibitem[Sa2]{sa:elliptic}
	K.~Saito, 
        {\it Extended affine root systems. I--IV}, 
	Publ. Res. Inst. Math. Sci. 21 (1985), %no. 1,
	75--179; 
	26 (1990)%, no. 1, 
	15--78;
	33 (1997)%, no. 2, 
 	301--329; 
	36 (2000) %, no. 3, 
	385--421. 
\bibitem[Sa3]{sa:duality}
	K.~Saito,
	{\it Duality for Regular Systems of Weights},
	Asian. J. Math. {\bf 2} no.4 (1998) 983-1048.
\bibitem[Sa4]{sa:around}
	K.~Saito,
	{\it Around the Theory of the Generalized Weight System: Relations 
	with Singularity Theory, the Generalized Weyl Group and 
	Its Invariant Theory, Etc.}, 
        (Japanese) Sugaku 38 (1986), %no. 2, no. 3, 
	97--115, 202--217; (Translation in English), 
	Amer. Math. Soc. Transl. (2) Vol.183 (1998) 101- 143.
\bibitem[Sa5]{sa:principal}
	K.~Saito, 
	{\it Principal $\Gamma$-cone for a tree $\Gamma$}, 
        preprint, RIMS-1507, June 2005, math.CO/0510623. 
\bibitem[Sc]{sch}
        F.~O.~Schreyer, 
        {\it Finite and countable CM-representation type}, 
        Singularities, representation of algebras, and vector bundles 
        (Lambrecht, 1985), 9--34, 
        Lecture Notes in Math., 1273, Springer, Berlin, 1987. 
\bibitem[T1]{t:1}
	A.~Takahashi,
	{\it K.~Saito's Duality for Regular Weight Systems and Duality for 
	Orbifoldized Poincar\'{e} Polynomials},
	Commun. Math. Phys. {\bf 205} (1999) 571-586.
\bibitem[T2]{t:2}
	A.~Takahashi,
	{\it Matrix Factorizations and Representations of Quivers I},
	math.AG/0506347. 
\bibitem[W]{w:1}
	J.~Walcher,
	{\it Stability of Landau-Ginzburg branes},
	J.Math.Phys. 46 (2005) 082305, hep-th/0412274. 
\bibitem[Y]{y}
        Y.~Yoshino, 
        {\it Cohen-Macaulay modules over Cohen-Macaulay rings}, 
        London Mathematical Society Lecture Note Series, 146, 
        Cambridge University Press, Cambridge, 1990. viii+177 pp. 











\end{thebibliography}
\end{document}